\documentclass[10pt]{article}
\usepackage{amsmath}
\usepackage{amsthm}
\usepackage{amssymb}
\usepackage{amstext}
\usepackage{epsfig}
\usepackage{epic}
\usepackage{algorithm}
\usepackage{algorithmic}

\newtheorem{theorem}{Theorem}[section]
\newtheorem{lemma}[theorem]{Lemma}

\newtheorem{definition}[theorem]{Definition}
\newtheorem{corollary}[theorem]{Corollary}

\newtheorem{example}[theorem]{Example}

\newcommand{\R}{\mathbb{R}}
\newcommand{\Q}{\mathbb{Q}}
\newcommand{\Z}{\mathbb{Z}}
\newcommand{\N}{\mathbb{N}}

\newcommand{\CP}{\mathcal{CP}}
\newcommand{\NF}{\mathcal{NF}}

\newcommand{\B}[2]{\mathcal{B}_{#1}(#2)}

\newcommand{\GRAPH}[3]{\mathcal{G}_{#1}(#2,#3)}
\newcommand{\Lattice}{\mathcal{L}}
\newcommand{\Feasible}{\mathcal{F}}
\newcommand{\FEASIBLE}[2]{\mathcal{F}_{#1}(#2)}
\newcommand{\FEASIBLEX}[3]{\mathcal{F}^{#1}_{#2}(#3)}
\newcommand{\IP}[3]{IP_{#1,#2}(#3)}
\newcommand{\IPX}[4]{IP^{#1}_{#2,#3}(#4)}
\newcommand{\seti}{{i}}
\newcommand{\inv}{\textnormal{-1}}
\newcommand{\zero}{{\bf 0}}

\newcommand{\sigmabar}{{\bar{\sigma}}}

\newcommand{\omegabar}{{\text{-}\omega}}
\newcommand{\ebar}{{\text{-}e}}
\newcommand{\f}{{\nu}}

\DeclareMathOperator{\rank}{rank}

\newcommand{\Kw}[1]{\textbf{#1}}

\begin{document}
\title{\bf Truncated Markov bases and Gr\"obner bases \\
for Integer Programming}
\author{
{\bf Peter N. Malkin} \\
\\\emph{CORE and INMA, Universit\'e catholique de Louvain, Belgium}\\
\texttt{malkin@core.ucl.ac.be}
}
\date{Version of \today}
\maketitle

\begin{abstract}
We present a new algorithm for computing a truncated Markov basis of a lattice.
In general, this new algorithm is faster than existing methods.
We then extend this new algorithm so that it solves the linear integer
feasibility problem with promising results for equality knapsack problems.
We also present a novel Gr\"obner basis approach to solve a particular integer
linear program as opposed to previous Gr\"obner basis methods that effectively
solved many different integer linear programs simultaneously.
Initial results indicate that this optimisation algorithm
performs better than previous Gr\"obner basis methods.


\end{abstract}


\section{Introduction}
Consider the set 
$\FEASIBLE{\Lattice}{\f}:=\{x\in\N^n: x\equiv \f\pmod{\Lattice}\}$
which we call a fiber of the lattice $\Lattice$
where $\f \in \Z^n$ and $\Lattice$ is a sub-lattice of $\Z^n$.
Importantly, the set of feasible solutions to any integer linear program
can be equivalently represented in the form $\FEASIBLE{\Lattice}{\f'}$ for some
lattice $\Lattice$ and some $\f' \in \Z^n$, and conversely, any fiber
$\FEASIBLE{\Lattice}{\f'}$ can be equivalently represented as the set
of feasible solutions to an integer linear program
(see Appendix \ref{Appendix: Lattice Programs and Integer Programs}).
We use the lattice notation of $\FEASIBLE{\Lattice}{\f}$ in this paper since we
find it notationally more convenient.

Informally, a Markov basis of the fiber $\FEASIBLE{\Lattice}{\f}$ is a
finite set of vectors in $\Lattice$ such that we can move from any feasible
solution to any other feasible solution in a finite number of steps via other
feasible solutions using the vectors in the Markov basis.
We can step from one feasible solution to another by adding or subtracting a
vector in the Markov basis.
A Markov basis of the lattice $\Lattice$ is a finite set of vectors
that is simultaneously a Markov basis for every possible fiber
$\FEASIBLE{\Lattice}{\f}$ as $\f$ varies over $\Z^n$.
Every lattice has a finite Markov basis, and so, every fiber has a finite
Markov basis.

One of the applications of Markov bases is in
algebraic statistics to test validity of statistical models via sampling
(see Diaconis and Sturmfels \cite{Diaconis+Sturmfels:98}).
Another area of application is in computational biology for problems arising
from phylogenetic trees (see \cite{Eriksson:04}).
Markov bases are also needed to perform integer optimisation using
Gr\"obner basis methods.
More and more problems are being solved with the growing computational power of
computer programs such as \texttt{4ti2},
an open source software package for algebraic, geometric, and
combinatorial problems on linear spaces  (\cite{4ti2}).

In some special situations, finding a Markov basis is straight-forward
(see \cite{Thomas+Weismantel:97}),
but in general, this is not the case.
There are three main methods for computing a Markov basis of a lattice,
(also called a generating set of a lattice):
the algorithm of Hosten and Sturmfels in \cite{Hosten+Sturmfels} called
the ``Saturation'' algorithm; the algorithm of Bigatti, LaScala, and
Robbiano in \cite{Bigatti+Lascala+Robbiano:99} that we call the
``Lift-and-Project'' algorithm, and the algorithm of Hemmecke and Malkin in
\cite{Hemmecke+Malkin:2006} called the ``Project-and-Lift'' algorithm.
Computationally, the Project-and-Lift algorithm is in general the fastest of
the three algorithms (\cite{Hemmecke+Malkin:2006}).

The above algorithms compute a Markov basis for all fibers of a lattice,
but we may only need a Markov basis for one fiber;
furthermore, there may be a huge difference between the size of a minimal
Markov basis of a lattice and the size of a minimal Markov basis of a fiber.
This makes it computationally worthwhile to focus on computing a Markov basis
for one fiber.
A truncated Markov basis with respect to the fiber $\FEASIBLE{\Lattice}{\f}$
is a Markov basis of a lattice (all fibers) after removing all vectors that
cannot be used to step between two feasible solutions in the fiber; in other
words, we remove vectors $u \in \Lattice$ when there do not exist 
$x,y \in \FEASIBLE{\Lattice}{\f}$ such that $x-y=u$.
A vector that cannot step between two feasible solutions in the fiber is never
needed in the Markov basis of a fiber, and thus,
a truncated Markov basis is thus a Markov basis of the fiber
$\FEASIBLE{\Lattice}{\f}$ but not every fiber.
A vector that can step between two feasible solutions in
$\FEASIBLE{\Lattice}{\f}$ still may not be
necessary in a Markov basis of $\FEASIBLE{\Lattice}{\f}$, 
but hopefully, there are not too many unnecessary vectors in the truncated
Markov basis.

The new algorithm for computing truncated Markov bases presented in this
paper in Section \ref{Section: Truncated Markov bases} is based upon the
Project-and-Lift algorithm combined with a truncated Gr\"obner basis algorithm
(see \cite{Thomas+Weismantel:97}), which is
described in Section \ref{Section: Computing Truncated Groebner bases}.
Previously, there were two known methods.
The first straight-forward method computes a Markov basis of a
lattice and then truncates the Markov basis.
The second method computes a truncated Graver basis of the
lattice (\cite{Hemmecke:Hilbert+Rays}), which is a superset of a truncated
Gr\"obner basis;
however, a truncated Graver basis might be a lot larger than a minimal
truncated Markov basis.
The new method is in general much faster than these two methods.

In Section \ref{Section:Feasibility}, we present a new algorithm to compute a
feasible solution of a fiber $\FEASIBLE{\Lattice}{\f}$
(i.e. find a point $x \in \FEASIBLE{\Lattice}{\f}$ if one exists) by using an
extension of the truncated Markov basis algorithm.
At the same time as finding a Markov basis, the truncated Markov basis 
algorithm can also compute a feasible solution without much additional
computation.
The notions of Markov bases and feasibility are indeed strongly related.
Recall that finding a feasible solution of a fiber is equivalent to finding a
feasible solution of a linear integer program.
Feasible solutions are not only interesting in themselves,
but are also needed to perform integer programming optimisation using a
Gr\"obner basis (see below).
We also show how to compute a feasible solution of a finite set of
different fibers $\FEASIBLE{\Lattice}{\f}$ simultaneously.
This is also extended to the truncated case where we focus on one
particular fiber.
We present the promising results of applying the algorithm to solve
equality constrained integer knapsacks (\cite{AardalLenstra2004}).

In Section \ref{Section:Optimality}, we solve linear integer programs using
Gr\"obner bases, which in this context are also known as \emph{test sets} for
integer linear programs.
Test sets were first introduced by Graver in \cite{Graver:75}
(see for example \cite{Weismantel:98}).
A Gr\"obner basis of a fiber $\FEASIBLE{\Lattice}{\f}$ with respect to some term
order $\succ$ is a finite set of vectors in $\Lattice$ such that
we can improve every non $\succ$-minimal point $x \in \FEASIBLE{\Lattice}{\f}$
by subtracting a vector in the Gr\"obner basis.
A Gr\"obner basis is thus a set of \emph{augmenting} or \emph{improving}
vectors.
So, using a $\succ$-Gr\"obner basis,
given some initial feasible solution of the fiber $\FEASIBLE{\Lattice}{\f}$,
we can find the $\succ$-minimal point in the fiber $\FEASIBLE{\Lattice}{\f}$
by iteratively improving the feasible solution until we can no longer do so,
in which case, we must have a $\succ$-minimal solution.
We call the problem of finding the $\succ$-minimal point in a fiber a
lattice program; that is, the problem
$\IP{\Lattice}{\succ}{\f} := \min_\succ\{x: x \in \FEASIBLE{\Lattice}{\f}\}$.
A Gr\"obner basis of a lattice is a finite set of vectors in $\Lattice$ that is
simultaneously a Gr\"obner basis for every fiber $\FEASIBLE{\Lattice}{\f}$ as
$\f$ varies over $\Z^n$.
Any integer linear program
can be expressed as a lattice program and conversely
any lattice program may be written as a integer linear program, and thus, the
concepts of integer linear programs and lattice programs are equivalent
(see Appendix \ref{Appendix: Lattice Programs and Integer Programs}).

The first Gr\"obner basis methods for solving a linear integer program 
constructed Gr\"obner bases of lattices and thus effectively solved
the program for all fibers simultaneously
(\cite{Conti+Traverso:91}).
This makes the approach appealing if we wish to solve the program for many
different fibers; however, if we wish to solve the program for only one fiber
then this is a disadvantage because we often do much more work than necessary.
By using the structure of a fiber, the new algorithm presented here focuses on
solving the program for a specific fiber.
To achieve this, the algorithm does two things: it solves a hierarchy of group
relaxations to avoid computing with unnecessary constraints
and it also applies truncation methods which are
strengthened by using the cost function.
Solving a hierarchy of group relaxations was explored in 
\cite{Gomory1965,Wolsey:71,Thomas+Hosten:03,Thomas:02}, and truncated methods
were explored in \cite{Thomas+Weismantel:97}; however, combining these two
approaches and strengthening truncation by using the cost function has not been
done before.

Gr\"obner bases and Markov bases have corresponding concepts in computational
algebraic geometry (see \cite{Conti+Traverso:91}) although we do not present it
here.
Instead, the approach in this paper follows the geometric approach in
\cite{Hemmecke+Malkin:2006,Thomas+Weismantel:97,Thomas:95,
Urbaniak+Weismantel+Ziegler:97, Weismantel:98}.

All computations in the paper were done using \texttt{4ti2} version 1.3
on a Intel Pentium4 3.0GHz machine running Linux. 
All timings given are rounded to the nearest one hundredth of a second.

\section{Truncated Markov Bases}
In this section, we define Markov bases of fibers, Markov bases of
lattices, and Markov bases of Gr\"obner bases.
This notion of truncation has been explored in
\cite{Urbaniak+Weismantel+Ziegler:97} and \cite{Thomas+Weismantel:97}, but only
for computing truncated Gr\"obner bases, and we apply it here to Markov bases.

A \textbf{lattice} is a set $\Lattice \subseteq \Z^n$ where
$\Lattice = \Z(S) = 
\{\sum^k_{i=1} \lambda_is^i : \lambda \in \Z^k\}$
for some finite set $S = \{s^1,\dots,s^k\} \subseteq \Z^n$.
If $\Lattice = \Z(S)$, then we say that $S$ spans $\Lattice$, and if $S$ is
inclusion-minimal, then we call $S$ a basis of $\Lattice$.

Given a lattice $\Lattice \subseteq \Z^n$, and a vector $\f \in \Z^n$,
we define the set
\[\FEASIBLE{\Lattice}{\f}:=\{x\in\N^n: x\equiv \f\pmod{\Lattice}\}
                       =\{x\in\N^n: x-\f \in \Lattice\}
\]
that we call a \textbf{fiber} of the lattice $\Lattice$.

Given a lattice $\Lattice \subseteq \Z^n$, a vector $\f \in \Z^n$, and a
set $S\subseteq\Lattice$, we define the \emph{fiber graph}
$\GRAPH{\Lattice}{\f}{S}$ to be the
\emph{undirected} graph with nodes $\FEASIBLE{\Lattice}{\f}$ and edges $(x,y)$
if $x-y\in S$ or $y-x\in S$ for $x,y \in \FEASIBLE{\Lattice}{\f}$.
\begin{definition}
Given $\f \in \Z^n$, we call a set $S\subseteq\Lattice$ a \emph{\bf Markov
basis of $\FEASIBLE{\Lattice}{\f}$} if the graph $\GRAPH{\Lattice}{\f}{S}$
is connected. The set $S$ is called a \emph{\bf Markov basis of $\Lattice$}
if it is a Markov basis of $\FEASIBLE{\Lattice}{\f}$ for every $\f \in \Z^n$.
\end{definition}
We remind the reader that connectedness of
$\GRAPH{\Lattice}{\f}{S}$ simply states that between
each pair $x,y\in \FEASIBLE{\Lattice}{\f}$ there exists a path
from $x$ to $y$ in the graph $\GRAPH{\Lattice}{\f}{S}$.

\begin{example}
Let $S := $\emph{\{(1,-1,-1,-3,1,2),(1,0,2,-2,-1,1)\}}, and
let $\Lattice \subseteq \Z^6$ be the lattice spanned by $S$.
By definition, $S$ is a spanning set of $\Lattice$, but $S$ is not a
Markov basis of $\Lattice$.
Observe that $\Lattice = \Lattice_A := \{u:Au=\zero, u\in\Z^n\}$ where
\[
A = (\tilde{A},I),\;
\tilde{A} = 
\begin{pmatrix}
 -2 & -3 \\
 +2 & -1 \\
 +1 & +2 \\
 -1 & +1  
\end{pmatrix}, \text{ and }
I = 
\begin{pmatrix}
 1 & 0 & 0 & 0 \\
 0 & 1 & 0 & 0 \\
 0 & 0 & 1 & 0 \\
 0 & 0 & 0 & 1  
\end{pmatrix}.
\]
So, for every $\f \in \Z^6$,
$\FEASIBLE{\Lattice}{\f} = F_A(b) =
\{(x,s): \tilde{A}x + Is = b, x\in\N^2,s\in\N^4\}$ where
$b = A\f \in \Z^4$.
Hence, the projection of $\FEASIBLE{\Lattice}{\f}$ onto the
$(x_1,x_2)$-plane is the set of integer points in the polyhedron
$\{x \in \R^n_+: \tilde{A}x \le b\}$,
and the $s$ variables are the slack variables.
Consider $\f := (2,2,4,2,5,1)$;
then, $\FEASIBLE{\Lattice}{\f} = F_A(b)$ where
$b = A\f = $\emph{(-6,4,11,1)}
(see Figure \ref{Figure: Connected Graph}a).

\begin{figure}[ht!]
\begin{center}
\framebox{ \setlength{\unitlength}{0.4pt}
\begin{picture}(310,260)(-20,-30)
    \put(260,210){(a)}
    \thicklines
    \put(0,0){\vector(0,1){220}}
    \multiputlist(0,-20)(50,0){0,1,2,3,4,5}
    \put(265,-20){$x_1$}
    \put(0,0){\vector(1,0){270}}
    \multiputlist(-20,0)(0,50){0,1,2,3,4}
    \put(-30,220){$x_2$}
    \matrixput(0,0)(50,0){6}(0,50){5}{\circle{7}}
    \dottedline{5}(0,50)(150,200)
    \dottedline{5}(0,100)(150,0)
    \dottedline{5}(100,0)(200,200)
    \dottedline{5}(150,200)(250,150)
    \put(100,100){\circle*{7}}
    \put(100,50){\circle*{7}}
    \put(50,100){\circle*{7}}
    \put(150,100){\circle*{7}}
    \put(100,150){\circle*{7}}
    \put(150,150){\circle*{7}}
    \put(150,200){\circle*{7}}
\end{picture}}
\framebox{ \setlength{\unitlength}{0.4pt}
\begin{picture}(310,260)(-20,-30)
    \put(260,210){(b)}
    \thicklines
    \put(0,0){\vector(0,1){220}}
    \multiputlist(0,-20)(50,0){0,1,2,3,4,5}
    \put(265,-20){$x_1$}
    \put(0,0){\vector(1,0){270}}
    \multiputlist(-20,0)(0,50){0,1,2,3,4}
    \put(-30,220){$x_2$}
    \matrixput(0,0)(50,0){6}(0,50){5}{\circle{7}}
    \dottedline{5}(0,50)(150,200)
    \dottedline{5}(0,100)(150,0)
    \dottedline{5}(100,0)(200,200)
    \dottedline{5}(150,200)(250,150)
    \put(100,100){\circle*{7}}
    \put(100,50){\circle*{7}}
    \put(50,100){\circle*{7}}
    \put(150,100){\circle*{7}}
    \put(100,150){\circle*{7}}
    \put(150,150){\circle*{7}}
    \put(150,200){\circle*{7}}

    \put(100,150){\line(1,0){50}}
    \put(100,100){\line(1,0){50}}
    \put(50,100){\line(1,0){50}}

    \put(150,100){\line(-1,1){50}}
    \put(100,50){\line(-1,1){50}}
\end{picture}}
\framebox{ \setlength{\unitlength}{0.4pt}
\begin{picture}(310,260)(-20,-30)
    \put(260,210){(c)}
    \thicklines
    \put(0,0){\vector(0,1){220}}
    \multiputlist(0,-20)(50,0){0,1,2,3,4,5}
    \put(265,-20){$x_1$}
    \put(0,0){\vector(1,0){270}}
    \multiputlist(-20,0)(0,50){0,1,2,3,4}
    \put(-30,220){$x_2$}
    \matrixput(0,0)(50,0){6}(0,50){5}{\circle{7}}
    \dottedline{5}(0,50)(150,200)
    \dottedline{5}(0,100)(150,0)
    \dottedline{5}(100,0)(200,200)
    \dottedline{5}(150,200)(250,150)
    \put(100,100){\circle*{7}}
    \put(100,50){\circle*{7}}
    \put(50,100){\circle*{7}}
    \put(150,100){\circle*{7}}
    \put(100,150){\circle*{7}}
    \put(150,150){\circle*{7}}
    \put(150,200){\circle*{7}}

    \put(100,150){\line(1,1){50}}
    \put(100,100){\line(1,1){50}}
    \put(50,100){\line(1,1){50}}
    \put(100,50){\line(1,1){50}}

    \put(100,150){\line(1,0){50}}
    \put(100,100){\line(1,0){50}}
    \put(50,100){\line(1,0){50}}

    \put(150,100){\line(-1,1){50}}
    \put(100,50){\line(-1,1){50}}
\end{picture}}
\end{center}
\caption{The set $\FEASIBLE{\Lattice}{\f}$ and the graphs
$\GRAPH{\Lattice}{\f}{S}$ and $\GRAPH{\Lattice}{\f}{S'}$ projected onto the
$(x_1,x_2)$-plane.}

\label{Figure: Connected Graph}
\end{figure}
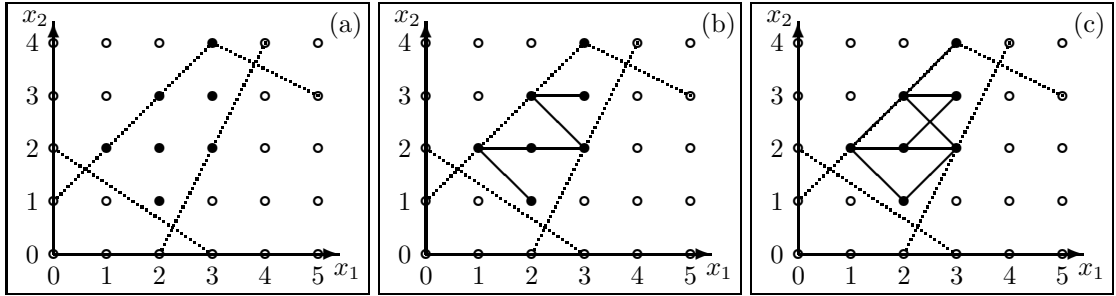

The graph
$\GRAPH{\Lattice}{\f}{S}$ is not connected
because the point $(3,4,12,2,0,0) \in \FEASIBLE{\Lattice}{\f}$ is disconnected
(see Figure \ref{Figure: Connected Graph}b). 
Let $S' := S \cup $\emph{\{(1,1,5,-1,-3,0)\}}. The graph of
$\GRAPH{\Lattice}{\f}{S'}$ is now connected
(see Figure \ref{Figure: Connected Graph}c).
Thus, $S'$ is a Markov basis of $\FEASIBLE{\Lattice}{\f}$.
\end{example}

A \textbf{truncated Markov basis} is a special type of Markov basis of
$\FEASIBLE{\Lattice}{\f}$ for some $\f \in \Z^n$ that is not necessarily a
Markov basis of all fibers.
Essentially, a truncated Markov basis with
respect to the fiber $\FEASIBLE{\Lattice}{\f}$ is a Markov basis of $\Lattice$
after removing all vectors $u \in \Lattice$ for which there does not exist $x,y
\in \FEASIBLE{\Lattice}{\f}$ such that $x-y=u$.
We call the act of removing such vectors \emph{truncation}.
More formally, let $G$ be a Markov basis of $\Lattice$; then, the set
$S := \{u\in G: u=x-y \text{ for some } x,y \in \FEASIBLE{\Lattice}{\f}\}$ 
is a truncated Markov basis.
Any vector $u \in \Lattice$ for which there does not exist
$x,y\in\FEASIBLE{\Lattice}{\f}$ such that $u=x-y$
can never be an edge in a fiber graph.
Hence, we never need such a vector $u$ in a Markov basis of
$\FEASIBLE{\Lattice}{\f}$.
Therefore, $S$ must be a Markov basis of the fiber $\FEASIBLE{\Lattice}{\f}$.

The set $S$ above is also a Markov basis of other related fibers.
Let $\f' \in \Z^n$ where $\FEASIBLE{\Lattice}{\f'} \ne \emptyset$
and $\FEASIBLE{\Lattice}{\f-\f'} \ne \emptyset$.
The set $S$ is also a Markov basis of the fiber
$\FEASIBLE{\Lattice}{\f'}$.
Let $u \in \Lattice$ for which there exists $x,y \in \FEASIBLE{\Lattice}{\f'}$
such that $u=x-y$, and let $\gamma \in \FEASIBLE{\Lattice}{\f-\f'}$.
Then, $x+\gamma,y+\gamma \in \FEASIBLE{\Lattice}{\f}$, and moreover,
$u = (x+\gamma)-(y+\gamma)$; thus, $u$ would not be removed during truncation.
So, any vector needed in a Markov basis of the fiber
$\FEASIBLE{\Lattice}{\f'}$ would not be removed by truncation, and therefore,
$S$ is still a Markov basis of $\FEASIBLE{\Lattice}{\f'}$.
The set $S$ is thus a Markov basis of the following set of fibers:
\[\B{\Lattice}{\f} :=  \{ \f' \in \Z^n: \FEASIBLE{\Lattice}{\f'} \ne \emptyset
\text{ and } \FEASIBLE{\Lattice}{\f-\f'} \ne \emptyset\}.\]
This property of a truncated Markov basis is the defining property
of truncated Markov bases.
\begin{definition}
Given $\f \in \Z^n$, we call a set $S\subseteq\Lattice$ a
\textbf{$\f$-truncated Markov basis of $\Lattice$} if 
$G$ is a Markov basis of $\FEASIBLE{\Lattice}{\f'}$ for every 
$\f' \in \B{\Lattice}{\f}$.
\end{definition}

Note that if $\FEASIBLE{\Lattice}{\f}\ne\emptyset$, then $\f \in \B{\Lattice}{\f}$
since $\FEASIBLE{\Lattice}{\zero}\ne\emptyset$.
Therefore, a $\f$-truncated Markov basis is by definition a Markov basis of
$\FEASIBLE{\Lattice}{\f}$, but a Markov basis of
$\FEASIBLE{\Lattice}{\f}$ is not necessarily a $\f$-truncated Markov basis.
Moreover, a $\f$-truncated Markov basis of $\Lattice$ is not necessarily a
Markov basis of $\Lattice$.
In the special case where $\FEASIBLE{\Lattice}{\f}=\emptyset$, we have
$\B{\Lattice}{\f} = \emptyset$, which is consistent since by definition an empty
set is a Markov basis of $\FEASIBLE{\Lattice}{\f}$ if
$\FEASIBLE{\Lattice}{\f}=\emptyset$.

Additionally, given a vector $u \in \Lattice$, there exists
$x,y\in\FEASIBLE{\Lattice}{\f}$ where $x-y=u$ if and only if
$\FEASIBLE{\Lattice}{\f-u^+} \ne \emptyset$
since $x-y = u$ means that
$x = \gamma+u^+$ and $y = \gamma+u^-$ for some $\gamma \in \N^n$ in which case
$\gamma \in \FEASIBLE{\Lattice}{\f-u^+} \ne \emptyset$.
Here, $u^+ \in \N^n$ is the positive part of $u$ and $u^-\in\N^n$ is the
negative part, or in other words,
$u^+_i = \max\{0,u_i\}$ and $u^-_i=\max\{0,-u_i\}$ for all
$i=1,...,n$.
Moreover, since $u^+ \in \FEASIBLE{\Lattice}{u^+} \ne \emptyset$,
we have that there exists $x,y\in\FEASIBLE{\Lattice}{\f}$
where $x-y=u$ if and only if $u^+\in\B{\Lattice}{\f}$.

The set $\B{\Lattice}{\f}$ has some interesting properties.
Given a set $S \subseteq \Lattice$, the connectivity of the graph
$\GRAPH{\Lattice}{\f'}{S}$ for the fibers in $\B{\Lattice}{\f}$ is strongly
related to the connectivity of $\GRAPH{\Lattice}{\f}{S}$.
Note that for $\f' \in \B{\Lattice}{\f}$,
we have $\gamma+\FEASIBLE{\Lattice}{\f'} \subseteq \FEASIBLE{\Lattice}{\f}$ for
every $\gamma \in \FEASIBLE{\Lattice}{\f-\f'}$ where
$\gamma+\FEASIBLE{\Lattice}{\f'} = \{\gamma+x:x \in \FEASIBLE{\Lattice}{\f'}\}$.
So, given $S \subseteq \Lattice$, any path $(x^0,...,x^k)$ in
$\GRAPH{\Lattice}{\f'}{S}$ can be translated by $\gamma$ to a path
$(x^0+\gamma,...,x^k+\gamma)$ in $\GRAPH{\Lattice}{\f}{S}$ for every
$\gamma \in \FEASIBLE{\Lattice}{\f-\f'}$.
Hence, if $S$ is a Markov basis of $\FEASIBLE{\Lattice}{\f'}$, then any
two points in $\gamma+\FEASIBLE{\Lattice}{\f'}\subseteq \FEASIBLE{\Lattice}{\f}$
are connected in $\GRAPH{\Lattice}{\f}{S}$.
Moreover, we have $\FEASIBLE{\Lattice}{\f'}+\FEASIBLE{\Lattice}{\f-\f'}\subseteq
\FEASIBLE{\Lattice}{\f}$ where 
$\FEASIBLE{\Lattice}{\f'}+\FEASIBLE{\Lattice}{\f-\f'}=
\{x+y: x \in \FEASIBLE{\Lattice}{\f'}, y \in \FEASIBLE{\Lattice}{\f-\f'}\}$.
Note that $\f-\f' \in \B{\Lattice}{\f}$ when $\f' \in \B{\Lattice}{\f}$.
If $S$ is both a Markov basis of $\FEASIBLE{\Lattice}{\f'}$ and a Markov basis
of $\FEASIBLE{\Lattice}{\f-\f'}$, then any two points in
$\FEASIBLE{\Lattice}{\f'}+\FEASIBLE{\Lattice}{\f-\f'}$ are 
connected in $\GRAPH{\Lattice}{\f}{S}$.
This is shown as follows.
Any two points in $\FEASIBLE{\Lattice}{\f'}+\FEASIBLE{\Lattice}{\f-\f'}$ can be
written in the form $x^1+y^1$ and $x^2+y^2$ where
$x^1,x^2 \in \FEASIBLE{\Lattice}{\f'}$ and
$y^1,y^2 \in \FEASIBLE{\Lattice}{\f-\f'}$.
Now, from above, the points $x^1+y^1, x^2+y^1 \in y^1+\FEASIBLE{\Lattice}{\f'}$
are connected in $\GRAPH{\Lattice}{\f}{S}$
and the points $x^2+y^1,x^2+y^2 \in x^2+\FEASIBLE{\Lattice}{\f-\f'}$
are connected in $\GRAPH{\Lattice}{\f}{S}$;
hence, the points $x^1+y^1$ and $x^2+y^2$ are connected in
$\GRAPH{\Lattice}{\f}{S}$ as required.

\begin{example}
Consider again the set
$S' :=$ \emph{\{(1,-1,-1,-3,1,2),(1,0,2,-2,-1,1),(1,1,5,-1,-3,0)\}}
and the lattice $\Lattice$ from above.
We saw previously that $S'$ is a Markov basis of
$\FEASIBLE{\Lattice}{\f}$ where $\f = (2,2,4,2,5,1)$;
however, $S'$ is not a $\f$-truncated Markov basis of $\Lattice$.
Consider $\f' = (2,2,4,2,0,0)$ (see Figure \ref{Figure: Disconnected Graph}a).
Note that $\f' \in \B{\Lattice}{\f}$ since
$\f-\f'=(0,0,0,0,5,1)\in\FEASIBLE{\Lattice}{\f-\f'} \ne \emptyset$.

\begin{figure}[ht!]
\begin{center}
\framebox{ \setlength{\unitlength}{0.4pt}
\begin{picture}(310,260)(-20,-30)
    \put(260,210){(a)}
    \thicklines
    \put(0,0){\vector(0,1){220}}
    \multiputlist(0,-20)(50,0){0,1,2,3,4,5}
    \put(265,-20){$x_1$}
    \put(0,0){\vector(1,0){270}}
    \multiputlist(-20,0)(0,50){0,1,2,3,4}
    \put(-30,220){$x_2$}
    \matrixput(0,0)(50,0){6}(0,50){5}{\circle{7}}
    \dottedline{5}(0,0)(200,200)
    \dottedline{5}(0,100)(150,0)
    \dottedline{5}(100,0)(200,200)
    \dottedline{5}(0,150)(250,25)
    \put(100,100){\circle*{7}}
    \put(100,50){\circle*{7}}
\end{picture}}
\framebox{ \setlength{\unitlength}{0.4pt}
\begin{picture}(310,260)(-20,-30)
    \put(260,210){(b)}
    \thicklines
    \put(0,0){\vector(0,1){220}}
    \multiputlist(0,-20)(50,0){0,1,2,3,4,5}
    \put(265,-20){$x_1$}
    \put(0,0){\vector(1,0){270}}
    \multiputlist(-20,0)(0,50){0,1,2,3,4}
    \put(-30,220){$x_2$}
    \matrixput(0,0)(50,0){6}(0,50){5}{\circle{7}}
    \dottedline{5}(0,50)(150,200)
    \dottedline{5}(0,50)(75,0)
    \dottedline{5}(100,0)(200,200)
    \dottedline{5}(0,50)(100,0)
    \put(0,50){\circle*{7}}
    \put(100,0){\circle*{7}}
\end{picture}}
\framebox{ \setlength{\unitlength}{0.4pt}
\begin{picture}(310,260)(-20,-30)
    \put(260,210){(c)}
    \thicklines
    \put(0,0){\vector(0,1){220}}
    \multiputlist(0,-20)(50,0){0,1,2,3,4,5}
    \put(265,-20){$x_1$}
    \put(0,0){\vector(1,0){270}}
    \multiputlist(-20,0)(0,50){0,1,2,3,4}
    \put(-30,220){$x_2$}
    \matrixput(0,0)(50,0){6}(0,50){5}{\circle{7}}
    \dottedline{5}(0,0)(200,200)
    \dottedline{5}(0,50)(75,0)
    \dottedline{5}(0,0)(100,200)
    \dottedline{5}(50,200)(250,100)
\end{picture}}

\end{center}
\caption{The sets $\FEASIBLE{\Lattice}{\f'}$, $\FEASIBLE{\Lattice}{\f''}$,
and $\FEASIBLE{\Lattice}{\f-\f''}$ projected onto the $(x_1,x_2)$-plane.}
\label{Figure: Disconnected Graph}
\end{figure}
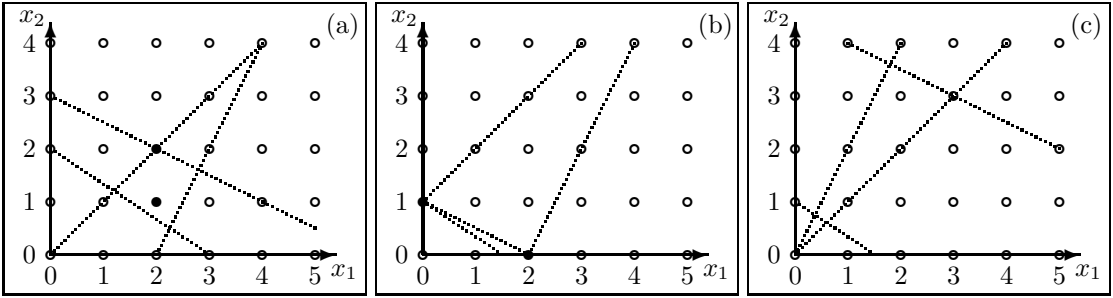

The graph $\GRAPH{\Lattice}{\f'}{S'}$ is disconnected since there are
only two feasible points in
$\FEASIBLE{\Lattice}{\f'} = \{(2,1,1,1,2,1),$ $(2,2,4,2,0,0)\}$.
The vector \emph{(0,1,3,1,-2,-1)} gives the unique minimal Markov basis of
$\FEASIBLE{\Lattice}{\f'}$, and hence, it must be in a $\f$-truncated
Markov basis of $\Lattice$.
The set $S'' = S' \cup $\emph{(0,1,3,1,-2,-1)} is a $\f$-truncated Markov basis
of $\Lattice$.

Although $S''$ is a $\f$-truncated Markov basis of $\Lattice$, it is not a
Markov basis for every $\f'' \in \Z^n$.
Consider $\f'' = $\emph{(0,1,0,5,0,0)}.
The graph $\GRAPH{\Lattice}{\f''}{S''}$ is disconnected since there are
only two feasible points in
$\FEASIBLE{\Lattice}{\f''} = \{(2,0,1,0,0,3),(0,1,0,5,0,0)\}$
(see Figure \ref{Figure: Disconnected Graph}b).
Observe that $\f'' \not\in \B{\Lattice}{\f}$ since
$\FEASIBLE{\Lattice}{\f-\f''} = \emptyset$
(see Figure \ref{Figure: Disconnected Graph}c).
\end{example}

In some situations, the $\f$-truncated Markov basis of $\Lattice$ may be
empty: a minimal $\f$-truncated Markov basis is empty if and only if
$|\FEASIBLE{\Lattice}{\f}| \le 1$ since then the fiber has zero or one elements.
On the other hand, a minimal $\f$-truncated Markov basis of $\Lattice$ may
also be a Markov basis of $\Lattice$.  So, the degree to which truncation
affects the size of the Markov basis varies between the two extremes of an
empty set and a Markov basis of a lattice.

\section{Truncated Gr\"obner bases}
\label{Section: Truncated Groebner bases}
In this section, we define Gr\"obner bases of fibers, Gr\"obner bases of
lattices, and truncated Gr\"obner bases
(see \cite{Urbaniak+Weismantel+Ziegler:97, Thomas+Weismantel:97}).

First, we need to define term orders.  We call $\succ$ a {\bf term ordering} for
$\Lattice$ if
\begin{enumerate}
\item $\succ$ is a total ordering on the set $\FEASIBLE{\Lattice}{\f}$ for
every $\f \in \Z^n$,
\item there is a unique $\succ$-minimal solution of
$\FEASIBLE{\Lattice}{\f}$ for every $\f \in \Z^n$ for which
$\FEASIBLE{\Lattice}{\f}\ne\emptyset$, and
\item $\succ$ is an additive ordering meaning that for all $\f \in \Z^n$ and for
all $x,y \in \FEASIBLE{\Lattice}{\f}$, if $x \succ y$, then $x+\gamma \succ
y+\gamma$ for every $\gamma\in\N^n$ (note that $x+\gamma, y+\gamma \in
\FEASIBLE{\Lattice}{\f+\gamma}$).
\end{enumerate}
The most common term orders are the lexicographic term ordering
and the degree reverse lexicographic term ordering
(see for example \cite{Cox+Little+OShea:92}). 

\begin{definition}
Given $\f \in \Z^n$ and a term order $\succ$, we call $G\subseteq\Lattice$ a
\emph{\bf $\succ$-Gr\"obner basis} of $\FEASIBLE{\Lattice}{\f}$
if for every $x\in\FEASIBLE{\Lattice}{\f}$, either $x$ is the unique
$\succ$-minimal element of $\FEASIBLE{\Lattice}{\f}$ or there exists a vector
$u \in G$ such that $x-u \in \FEASIBLE{\Lattice}{\f}$ and $x \succ x-u$.
The set $G$ is called a $\succ$-Gr\"obner basis of $\Lattice$ if it is a
Gr\"obner basis for every $\f\in\Z^n$
\end{definition}

Analogously to truncated Markov bases, we define truncated Gr\"obner bases.
\begin{definition}
Given $\f \in \Z^n$, we call a set $G\subseteq\Lattice$ a
\emph{\bf $\f$-truncated $\succ$-Gr\"obner basis of $\Lattice$} if 
$G$ is a Gr\"obner basis of $\FEASIBLE{\Lattice}{\f'}$ for every 
$\f' \in \B{\Lattice}{\f}$.
\end{definition}

As in the Markov basis case, a $\f$-truncated $\succ$-Gr\"obner basis is by
definition a $\succ$-Gr\"obner basis of
$\FEASIBLE{\Lattice}{\f}$, but a $\succ$-Gr\"obner basis of
$\FEASIBLE{\Lattice}{\f}$ is not necessarily a $\f$-truncated $\succ$-Gr\"obner
basis and furthermore, a $\f$-truncated $\succ$-Gr\"obner basis of $\Lattice$ is
not necessarily a $\succ$-Gr\"obner basis of $\Lattice$.

We can solve \emph{lattice programs} using Gr\"obner bases.
Given a lattice $\Lattice$, a vector $\f$, and a term order $\succ$, the
problem
\[\IP{\Lattice}{\succ}{\f}:=\min_\succ\{x:x\in\FEASIBLE{\Lattice}{\f}\}\]
is called a \textbf{lattice program}.
Given some initial feasible solution $x\in\FEASIBLE{\Lattice}{\f}$,
we can solve the lattice program $\IP{\Lattice}{\succ}{\f}$
using a Gr\"obner basis $G$ of the fiber $\FEASIBLE{\Lattice}{\f}$
by iteratively improving the feasible solution using vectors in $G$.
This process constructs a $\succ$-decreasing path in the graph
$\GRAPH{\Lattice}{\f}{G}$ from the initial feasible solution $x$ to the unique
$\succ$-minimal solution
where a path $(x^0,\ldots,x^k)$ in $\GRAPH{\Lattice}{\f}{G}$ is
{\bf $\succ$-decreasing} if $x^i\succ x^{i+1}$ for $i=0,\ldots,k-1$.
This gives us an equivalent way of defining Gr\"obner bases in terms of paths in
the graph $\GRAPH{\Lattice}{\f}{G}$ in Lemma
\ref{Lemma:decreasing paths} below.
Firstly, note that for a vector $u \in \Lattice$,
where $x,x-u \in \FEASIBLE{\Lattice}{\f}$ and $x \succ x-u$, we must have
$u^+ \succ u^-$.
This property follows since $\succ$ is an additive ordering.
Thus, we only need consider vectors in the set
$\Lattice_\succ := \{u \in \Lattice: u^+ \succ u^-\}$.
\begin{lemma}
\label{Lemma:decreasing paths}
Given $\f \in \Z^n$, $G\subseteq\Lattice_\succ$ is a
{\bf $\succ$-Gr\"obner basis} of $\FEASIBLE{\Lattice}{\f}$
if and only if for every $x\in\FEASIBLE{\Lattice}{\f}$ there exists a decreasing
path in $\GRAPH{\Lattice}{\f}{G}$
from $x$ to the unique $\succ$-minimal element in $\FEASIBLE{\Lattice}{\f}$.
\end{lemma}

Importantly, if $G\subseteq\Lattice_\succ$ is a $\succ$-Gr\"obner basis, then
$G$ is a Markov basis of $\Lattice$ since given $x,y \in
\GRAPH{\Lattice}{\f}{G}$ for some $\f \in \Z^n$, there exists a
$\succ$-decreasing path from $x$ to the unique
$\succ$-minimal element in $\FEASIBLE{\Lattice}{\f}$ and from $y$ to the same
element, and thus, $x$ and $y$ are connected in
$\GRAPH{\Lattice}{\f}{G}$.

We can use Gr\"obner bases to solve the 
integer program
\[\IP{\Lattice}{c}{\f}:=\min\{cx:x\in\FEASIBLE{\Lattice}{\f}\}\]
given a lattice $\Lattice$, a vector $\f \in \Z^n$,
and a cost vector $c \in \Z^n$.
To solve $\IP{\Lattice}{c}{\f}$, we solve instead a lattice program
$\IP{\Lattice}{\succ}{\f}$ for some \emph{term order} $\succ$ that is compatible
with $c$.

Given a vector $c \in \Z^n$, we say that a vector a term order $\succ$ is
\textbf{compatible} with $c$ if
the optimal solution of $\IP{\Lattice}{\succ}{\f}$ is also an optimal solution
of $\IP{\Lattice}{c}{\f}$ for all $\f\in\Z^n$ where
$\FEASIBLE{\Lattice}{\f}\ne\emptyset$.
We can easily construct a $c$ compatible order $\succ_c$ given some 
(tie-breaking) term ordering $\succ$ as follows: $x \succ_c y$ if $cx > cy$, or
$cx=cy$ and $x \succ y$.
We must be a little careful here though since $\succ_c$ is not necessarily a
term order.
The ordering $\succ_c$ satisfies conditions (i) and (iii) for being a term
order, but condition (ii) is not always satisfied.
The ordering $\succ_c$ is a term ordering if and only if
$\IP{\Lattice}{c}{\f}$ has an optimal solution
for every $\f \in \Z^n$ where $\FEASIBLE{\Lattice}{\f} \ne \emptyset$.
Note that $\IP{\Lattice}{c}{\f}$ has an optimal solution for every
$\f \in \Z^n$ where $\FEASIBLE{\Lattice}{\f} \ne \emptyset$ if and only if
$\IP{\Lattice}{c}{\zero} := \min\{cx:x\in\Lattice,x\in\N^n\} = 0$,
and we can check whether $\IP{\Lattice}{c}{\zero}=0$ using linear programming.


Conversely, given any term order $\succ$ there always exists a compatible
$c\in\Z^n$ such that $\IP{\Lattice}{c}{\f}$ has a unique optimal solution that
is also the optimal solution of $\IP{\Lattice}{\succ}{\f}$ for all $\f\in\Z^n$
where $\FEASIBLE{\Lattice}{\f}\ne\emptyset$ (\cite{Sturmfels:96}).
So, solving integer programs $\IP{\Lattice}{c}{\f}$ is essentially equivalent
to solving lattice programs.

\section{Computing truncated Gr\"obner bases}
\label{Section: Computing Truncated Groebner bases}
In this section, we describe how to compute truncated Gr\"obner bases.
We first describe how to compute truncated Markov bases since
the algorithm for computing truncated Gr\"obner bases is used to compute
truncated Markov bases.
We present existing algorithms for truncated Gr\"obner bases including a new 
approaches to truncation.
The structure of this section follows closely from 
\cite{Hemmecke+Malkin:2006} except that we now deal with truncation,
and so, we have omitted the proofs of those results
in this section that correspond closely to results in
\cite{Hemmecke+Malkin:2006}.

\begin{figure}[ht]
\begin{center}
\framebox{ \setlength{\unitlength}{0.4pt}
\begin{picture}(600,120)(-10,-10)

    \put(50,80){\circle*{5}}
    \put(100,10){\circle*{5}}
    \put(150,60){\circle*{5}}
    \put(250,20){\circle*{5}}
    \put(350,40){\circle*{5}}
    \put(400,50){\circle*{5}}
    \put(450,0){\circle*{5}}
    \put(550,30){\circle*{5}}

    \drawline(550,30)(450,0)(400,50)(350,40)(250,20)(150,60)(100,10)
        (50,80)

    \put(43,90){$x$}
    \put(543,40){$y$}

    \put(95,25){$x^1$}
    \put(140,65){$x^2$}
    \put(240,30){$x^3$}
    \put(343,50){$x^4$}
    \put(393,60){$x^5$}
    \put(443,15){$x^6$}

    \thicklines
    \put(0,0){\vector(0,1){90}}
    \put(-10,100){$\prec$}
\end{picture}}
\end{center}
\caption{Reduction path between $x$ and $y$.}

\label{Figure: Reduction Path}
\end{figure}
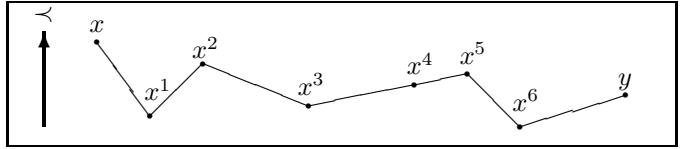

We now describe Gr\"obner bases in terms of reduction paths, so that we avoid
explicitly mentioning the $\succ$-minimal solution of every fiber.
A path $(x^0,\ldots,x^k)$ in $\GRAPH{\Lattice}{\f}{G}$
is a {\bf $\succ$-reduction path} if for all
$i\in\{1,\ldots,k-1\}$, we have either $x^0\succ x^i$ or $x^k\succ x^i$.
For example, see Figure \ref{Figure: Reduction Path}.

\begin{lemma} \label{Lemma: reduction paths}
Given $\f \in \Z^n$, a set $G\subseteq\Lattice_\succ$
is a $\succ$-Gr\"obner basis of $\FEASIBLE{\Lattice}{\f}$ if
and only if for each pair $x,y\in \FEASIBLE{\Lattice}{\f}$, there exists a
$\succ$-reduction path in $\GRAPH{\Lattice}{\f}{G}$ between $x$ and
$y$.
\end{lemma}
%

Checking for a given $G\subseteq\Lattice_\succ$ whether there
exists a $\succ$-reduction path in $\GRAPH{\Lattice}{\f}{G}$
for every $\f\in\Z^n$ and for each pair $x,y\in \FEASIBLE{\Lattice}{\f}$ involves
many situations that need to be checked. In fact, far fewer checks
are needed: we only need to check for a $\succ$-reduction path from $x$ to $y$
if there exists a \textbf{$\succ$-critical path} from $x$ to $y$.

\begin{figure}[ht!]
\begin{center}
\framebox{ \setlength{\unitlength}{0.2pt}
\begin{picture}(320,210)(-20,-20)

    \put(50,40){\circle*{10}}
    \put(150,140){\circle*{10}}
    \put(270,20){\circle*{10}}

    \drawline(50,40)(150,140)(270,20)

    \put(45,5){$x$}
    \put(135,155){$z$}
    \put(265,-15){$y$}

    \put(70,100){$u$}
    \put(220,85){$v$}

    \thicklines

    \put(0,0){\vector(0,1){160}}
    \put(-10,165){$\prec$}
\end{picture}}
\end{center}
\caption{A critical path for $(u,v)$ between $x$, $z$, and $y$.}

\label{Figure: uv Path}
\end{figure}
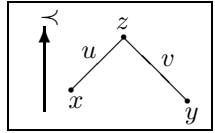

\begin{definition}
Given $G \subseteq \Lattice_\succ$ and $\f \in \Z^n$, a path $(x,z,y)$
in $\GRAPH{\Lattice}{\f}{G}$
is a {\bf $\succ$-critical path} if $z\succ x$ and $z\succ y$.
\end{definition}
If $(x,z,y)$ is a $\succ$-critical path in $\GRAPH{\Lattice}{\f}{G}$,
then $x+u=z=y+v$ for some pair $u,v \in G$, in which case, we call $(x,z,y)$ a
$\succ$-critical path for $(u,v)$ (see Figure \ref{Figure: uv Path}).

\begin{lemma}
\label{Lemma: critical reduction paths}
Given $\f \in \Z^n$, let $G\subseteq\Lattice_\succ$ where $G$ is
a Markov basis of $\FEASIBLE{\Lattice}{\f}$.
$G$ is a $\succ$-Gr\"obner basis of $\FEASIBLE{\Lattice}{\f}$ if and only if
there exists a $\succ$-reduction path between $x'$ and $y'$ for every
$\succ$-critical path $(x',z',y')$ in $\GRAPH{\Lattice}{\f}{G}$.
\end{lemma}

We can extend Lemma \ref{Lemma: critical reduction paths} to
$\f$-truncated Gr\"obner bases. It is a straight-forward consequence of
Lemma \ref{Lemma: critical reduction paths}, but nevertheless, it is worthwhile
stating explicitly.
\begin{lemma} \label{Lemma: truncated critical reduction paths}
Given $\f \in \Z^n$, let $G\subseteq\Lattice_\succ$ where $G$ is
a $\f$-truncated Markov basis of $\Lattice$.
$G$ is a $\f$-truncated $\succ$-Gr\"obner basis of $\Lattice$ if
and only if there exists a $\succ$-reduction path between $x'$ and $y'$ for
every $\succ$-critical path $(x',z',y')$ in $\GRAPH{\Lattice}{\f'}{G}$
for all $\f' \in \B{\Lattice}{\f}$.
\end{lemma}

It is not necessary to check for a $\succ$-reduction path from
$x'$ to $y'$ for every $\succ$-critical path $(x',y',z')$
in $\GRAPH{\Lattice}{\f'}{G}$
for all $\f'\in\B{\Lattice}{\f}$. Consider the case where there exists another
$\succ$-critical path $(x'',y'',z'')$ in $\GRAPH{\Lattice}{\f''}{G}$
for some $\f''\in\Z^n$ such that $(x',y',z')=(x''+\gamma,y''+\gamma,z''+\gamma)$
for some $\gamma \in \N^n$. Then, a $\succ$-reduction path from
$x''$ to $y''$ in $\GRAPH{\Lattice}{\f''}{G}$ translates by $\gamma$
to a $\succ$-reduction path from $x'$ to $y'$ in
$\GRAPH{\Lattice}{\f'}{G}$.
Moreover, $\f'' \in \B{\Lattice}{\f}$, since
$\gamma \in \FEASIBLE{\Lattice}{\f'-\f''}\ne\emptyset$ which implies that
$\FEASIBLE{\Lattice}{\f-\f''}\ne\emptyset$ because
$\FEASIBLE{\Lattice}{\f-\f'}\ne\emptyset$.
Thus, we only need to check for a $\succ$-reduction path from $x''$ to $y''$.

We call a $\succ$-critical path $(x,y,z)$ {\bf minimal} if there does not exist
another $\succ$-critical path $(x',y',z')$ such that
$(x,y,z)=(x'+\gamma,y'+\gamma,z'+\gamma)$
for some $\gamma \in \N^n$ where $\gamma \ne 0$, or equivalently,
$\min\{x_i,y_i,z_i\}=0$ for all $i = 1,\dots,n$.
Consequently, if there exists a $\succ$-reduction path between $x'$ and $y'$ for
all minimal $\succ$-critical paths $(x',y',z')$ in
$\GRAPH{\Lattice}{\f'}{G}$ for
some $\f' \in \B{\Lattice}{\f}$, then there exists a
$\succ$-reduction path between $x''$ and $y''$ for all $\succ$-critical paths
$(x'',y'',z'')$ in $\GRAPH{\Lattice}{\f''}{G}$ for
some $\f'' \in \B{\Lattice}{\f}$.
Also, for each pair of vectors $u,v \in G$, there
exists a unique minimal $\succ$-critical path $(x^{(u,v)},z^{(u,v)},y^{(u,v)})$
in $\GRAPH{\Lattice}{\f^{(u,v)}}{G}$
determined by $z^{(u,v)}:=\max\{u^+,v^+\}$ component-wise,
$x^{(u,v)}:=z^{(u,v)}-u$, $y^{(u,v)}:=z^{(u,v)}-v$, and
$\f^{(u,v)}:=z^{(u,v)}$.
So, any other $\succ$-critical path for $(u,v)$ is of the form
$(x^{(u,v)}+\gamma,z^{(u,v)}+\gamma,y^{(u,v)}+\gamma)$ for some
$\gamma\in\N^n$. Using minimal $\succ$-critical paths, we can rewrite
Lemma \ref{Lemma: truncated critical reduction paths}, so that we only need to
check for a \emph{finite} number of $\succ$-reduction paths.

\begin{lemma} \label{Lemma: minimal critical paths}
Let $G\subseteq\Lattice_\succ$ and $\f \in \Z^n$.
The set $G$ is a $\f$-truncated $\succ$-Gr\"obner basis of $\Lattice$ if and
only if for all $u,v \in G$ where $\f^{(u,v)} \in \B{\Lattice}{\f}$,
there exists a $\succ$-reduction path between $x^{(u,v)}$ and $y^{(u,v)}$
in $\GRAPH{\Lattice}{\f^{(u,v)}}{G}$.
\end{lemma}

We now turn Lemma \ref{Lemma: minimal critical paths} into an algorithmic tool.
Algorithm \ref{Algorithm: Completion procedure} below,
called a completion procedure (\cite{Buchberger:87}),
starts from a $\f$-truncated Markov basis and 
computes a $\f$-truncated $\succ$-Gr\"obner basis. 
An important part of this algorithm is checking whether
$\f' \in \B{\Lattice}{\f}$ for some $\f' \in \Z^n$.  How exactly we perform this
check in practice is discussed at length after first presenting the
overall algorithm.

Given a set $S \subseteq \Lattice$, the completion procedure first sets $G:=S$,
and then directs all vectors in $G$ according to $\succ$ such that
$G \subseteq \Lattice_\succ$. It also removes from the set $G$ any vectors
$u \in G$ such that $u^+ \not \in \B{\Lattice}{\f}$ -- recall that these vectors
are not needed in a truncated Markov basis.
Note that at this point
$\GRAPH{\Lattice}{\f'}{S} = \GRAPH{\Lattice}{\f'}{G}$
for all $\f' \in \B{\Lattice}{\f}$, and thus, $G$ is also a $\f$-truncated
Markov basis of $\Lattice$.
The completion procedure then determines whether the set $G$ satisfies Lemma
\ref{Lemma: minimal critical paths}; in other words, it tries to find a
reduction path from $x^{(u,v)}$ to $y^{(u,v)}$ for every pair $u,v \in G$
where $\f^{(u,v)} \in \B{\Lattice}{\f}$.
If $G$ satisfies Lemma \ref{Lemma: minimal critical paths}, then we are done.
Otherwise, no $\succ$-reduction path was found for some $(u,v)$, in which case,
we add a vector to $G$ so that a $\succ$-reduction path exists, and then
again, test whether $G$ satisfies Lemma
\ref{Lemma: minimal critical paths} and so on.

To check for a $\succ$-reduction path from $x^{(u,v)}$ to $y^{(u,v)}$,
we construct a maximal $\succ$-decreasing path
in $\GRAPH{\Lattice}{\f^{(u,v)}}{G}$ from
$x^{(u,v)}$ to some $x'$ and
from $y^{(u,v)}$ to some $y'$
using the ``Normal Form Algorithm''
(Algorithm \ref{Algorithm: normal form} below).
If $x' = y'$,  then we have found a $\succ$-reduction path from
$x^{(u,v)}$ to $y^{(u,v)}$. Otherwise, we add the vector $r \in \Lattice_\succ$
to $G$ where $r := x'-y'$ if $x' \succ y'$, and $r := y'-x'$ otherwise,
and then, there is now a $\succ$-reduction path from
$x^{(u,v)}$ to $y^{(u,v)}$ in $\GRAPH{\Lattice}{\f^{(u,v)}}{G}$.

\begin{algorithm}
\caption{Normal Form Algorithm}
\begin{algorithmic}
\REQUIRE a vector $x\in\N^n$ and a set $G\subseteq\Lattice_\succ$. 
\ENSURE a vector $x'$ where there is a maximal $\succ$-decreasing path
from $x$ to $x'$ in $\GRAPH{\Lattice}{x}{G}$.
\STATE $x' := x$
\WHILE{there is exists $u\in G$ such that $u^+\leq x'$}
  \STATE $x':=x'-u$
\ENDWHILE
\STATE \Kw{return} $x'$
\end{algorithmic}
\label{Algorithm: normal form}
\end{algorithm}

\begin{algorithm} 
\caption{Truncated completion procedure}
\begin{algorithmic}
\REQUIRE a vector $\f \in \Z^n$, a term ordering $\succ$,
and a $\f$-truncated Markov basis $S\subseteq\Lattice$.
\ENSURE a $\f$-truncated Gr\"obner basis $G\subseteq\Lattice_\succ$.
\STATE $G:=\{u:u^+\succ u^-, u\in S\}\cup\{-u:u^-\succ u^+, u\in S\}$
\STATE $G:=\{u:u\in G, u^+ \in \B{\Lattice}{\f}\}$
\STATE $C:=\{(u,v): u,v \in G, \f^{(u,v)} \in \B{\Lattice}{\f}\}$
\WHILE{$C\neq \emptyset$}
  \STATE Select $(u,v) \in C$
  \STATE $C:=C\setminus\{(u,v)\}$
  \STATE $r:=\NF(x^{(u,v)},G)-\NF(y^{(u,v)},G)$
  \IF{$r\neq 0$}
     \STATE \Kw{if} $r^- \succ r^+$ \Kw{then} $r:=-r$ 
     \STATE $C:=C\cup\{(r,s):s\in G, \f^{(r,s)}\in\B{\Lattice}{\f}\}$
     \STATE $G:=G\cup \{r\}$
  \ENDIF
\ENDWHILE
\STATE \Kw{return} $G$.
\end{algorithmic}
\label{Algorithm: Completion procedure}
\end{algorithm}
We write $\NF(x,G)$ for the output of the Normal Form Algorithm and 
we write $\CP_\f(\succ,S)$ for the output of the Completion Procedure.

%

There is a trade off between the computational benefit of computing a
$\f$-truncated Markov basis (computing a smaller set) and the computational
cost of computing whether $\f' \in \B{\Lattice}{\f}$ for some $\f' \in \Z^n$ many times.
In general, it is NP-hard to determine whether $\f' \in \B{\Lattice}{\f}$ since we must
know if $\FEASIBLE{\Lattice}{\f'}\ne \emptyset$ and
$\FEASIBLE{\Lattice}{\f-\f'}\ne \emptyset$.
Instead, we can check a sufficient condition for when
$\f' \not \in \B{\Lattice}{\f}$, and so, we compute a superset of a
$\f$-truncated Gr\"obner basis since we keep some vectors that are not needed.

Firstly, note that in Algorithm \ref{Algorithm: Completion procedure}, whenever
we check whether $\f' \in \B{\Lattice}{\f}$, we always have
$\FEASIBLE{\Lattice}{\f'}\ne\emptyset$ since either $\f'=u^+\ge 0$ for some
$u \in \Lattice$ or $\f'=\f^{(u,v)} \ge 0$ for some
$u,v \in \Lattice$ and in either
case $\f' \in \FEASIBLE{\Lattice}{\f'} \ne \emptyset$. So, in the algorithm, we
only need to check whether $\FEASIBLE{\Lattice}{\f-\f'}\ne \emptyset$.

We could instead check for feasibility of a relaxation of the feasible set
$\FEASIBLE{\Lattice}{\f-\f'}$.
One possible relaxation of $\FEASIBLE{\Lattice}{\f-\f'}$ to check is 
\[\FEASIBLEX{\Z}{\Lattice}{\f-\f'} :=\{x: x\equiv \f-\f'\pmod{\Lattice},x\in\Z^n\}.\] 
But since $\f-\f' \in \Z^n$,
we have $\f-\f' \in \FEASIBLEX{\Z}{\Lattice}{\f-\f'} \ne \emptyset$, and so, this is
trivially always satisfied.
Another possible relaxation is the linear
relaxation of $\FEASIBLE{\Lattice}{\f-\f'}$:
\[\FEASIBLEX{\R}{\Lattice}{\f-\f'}:=\{x:x\equiv \f-\f'\pmod{\Lattice^\R},x\in\R^n_+\}\]
where $\Lattice^\R\subseteq\R^n$ is the smallest
subspace containing $\Lattice$; that is,
$\Lattice^R := \{ku:u \in \Lattice, k \in \R\}$.
We can thus solve a linear program to check whether
$\FEASIBLEX{\R}{\Lattice}{\f-\f'} = \emptyset$ implying that
$\FEASIBLE{\Lattice}{\f-\f'} = \emptyset$.
Note that $\FEASIBLE{\Lattice}{\f-\f'} = \FEASIBLEX{\R}{\Lattice}{\f-\f'} \cap
\FEASIBLEX{\Z}{\Lattice}{\f-\f'}$, but
 $\FEASIBLEX{\R}{\Lattice}{\f-\f'} \ne \emptyset$ and
$\FEASIBLEX{\Z}{\Lattice}{\f-\f'} \ne \emptyset$ do not imply that
$\FEASIBLE{\Lattice}{\f-\f'} \ne \emptyset$.

In practice, computational experiments
show that it is usually not worthwhile performing the full check
whether $\FEASIBLEX{\R}{\Lattice}{\f-\f'} = \emptyset$, so instead,
we use a sufficient condition for when
$\FEASIBLEX{\R}{\Lattice}{\f-\f'} = \emptyset$ and thus $\f' \not \in
\B{\Lattice}{\f}$ that is quick to check.
Let \[\Lattice^\circ_+ := \{a \in \R^n_+: au = 0 \; \forall u \in \Lattice\}.\]
Note that $\Lattice^\circ_+$ is a pointed convex cone.
Firstly, observe that, for any $a \in \Lattice^\circ_+$ and any $\f' \in \Z^n$, we
have $ax = a\f'$ for all $x \in \FEASIBLE{\Lattice}{\f'}$.
Secondly, if $a \in \Lattice^\circ_+$,
then, for all $\f' \in \B{\Lattice}{\f}$, we have $a\f \ge a\f'$
since if $x \in \FEASIBLE{\Lattice}{\f-\f'}\ne \emptyset$, then $a(\f-\f') = ax$
and $ax \ge 0$ because $a \ge 0$ and $x \ge 0$. 
Therefore, if $a\f < a\f'$, then $\FEASIBLEX{\R}{\Lattice}{\f-\f'} = \emptyset$
and thus $\f' \not \in \B{\Lattice}{\f}$.
Moreover, it follows from Farkas' lemma that
$\FEASIBLEX{\R}{\Lattice}{\f-\f'} = \emptyset$ if and only if there exists an
$a \in \Lattice^\circ_+$ where $a\f < a\f'$, and furthermore,
$\FEASIBLEX{\R}{\Lattice}{\f-\f'} = \emptyset$ if and only if there exists an
extreme ray $a$ of the cone $\Lattice^\circ_+$ where $a\f < a\f'$.
The set of extreme rays is finite but there are far too many of them in general
to check this condition. So, we need a way of selecting one
 $a \in \Lattice^\circ_+$ or a small set of $a$.
Choosing different $a$'s can produce very different results, and
the best $a$'s to choose vary from fiber to fiber.

We now present a novel approach for selecting a \emph{good}
$a \in \Lattice^\circ_+$.
Now, note that when we run Algorithm \ref{Algorithm: Completion procedure} and
check whether $\f' \in \B{\Lattice}{\f}$, we have $a\f' \ge 0$ for all
$a \in \Lattice^\circ_+$ since from above 
$\FEASIBLE{\Lattice}{\f'} \ne \emptyset$ ($\f' \ge 0$).
Ideally, there exists $a \in \Lattice^\circ_+$ where $a\f = 0$,
implying that $a\f'=0$ (i.e. if $a_i \ne 0$, then $\f'_i = 0$) for every
$\f' \in \B{\Lattice}{\f}$.
This condition is very strong and is quick to check and effectively means that
we compute using a sub-lattice of $\Lattice$.
Otherwise if $a\f > 0$ for all $a \in \Lattice^\circ_+$, then
a useful heuristic is to choose a single $a \in \Lattice^\circ_+$ such
that $a\f$ is minimal with respect to some norm $||\cdot||$ of $a$. More
formally, we solve the following problem:
\[ \text{argmin} \{ a\f : ||a|| = 1, a \in \Lattice^\circ_+\}.\]
If we use the $l_1$-norm (i.e. $||a||_1 = \sum_i a_i$), then we can find $a$
using linear programming.
In this case, note that we only need to solve one linear program to compute $a$
as opposed to solving a linear program every time we check whether
$\f'\in \B{\Lattice}{\f}$.

\begin{example}
Consider again the set $S := $\emph{\{(1,-1,-1,-3,1,2),(1,0,2,-2,-1,1)\}} and
the lattice $\Lattice \subseteq \Z^6$ spanned by $S$.
Let $c:=(2,1,0,0,0,0)$ and let $\prec$ be some term order.
Then, a $\prec_c$-Gr\"obner basis of $\Lattice$ is
\[G := \text{\emph{\{(1,-1,-1,-3,1,2), (1,0,2,-2,-1,1), (1,1,5,-1,-3,0),
(0,1,3,1,-2,-1), (2,-1,1,-5,0,3)\}}}.\]
Recall that this is a $\prec_c$-Gr\"obner basis for every possible fiber
$\f\in\Z^n$.
We now examine truncated $\prec_c$-Gr\"obner basis for two different fibers.

\begin{enumerate}
\item
Consider $\f:=(0,1,0,5,0,0)$ (see Figure \ref{Figure: Disconnected Graph}b);
$\FEASIBLE{\Lattice}{\f}=\{(2,0,1,0,0,3),(0,1,0,5,0,0)\}$.
Since the feasible set consists of only two feasible solutions, the minimal
$\f$-truncated $\prec_c$-Gr\"obner basis contains only one vector:
the vector between the two feasible solutions.
Thus, the set \emph{\{(2,-1,1,-5,0,3)}\} is a $\f$-truncated $\prec_c$-Gr\"obner
basis of $\Lattice$.

If we run the truncated completion procedure,
Algorithm \ref{Algorithm: Completion procedure}, using
$\FEASIBLEX{\R}{\Lattice}{\f-\f'} \ne \emptyset$ as a check for truncation, we
compute the set $\CP_\f(\succ,S) = \emph{\{(2,-1,1,-5,0,3)}\}$.

If instead we run the truncated completion procedure using the quick truncation
test $a\f < a\f'$ where
$a := \frac{1}{4}(0,1,1,0,2,0)=
\text{argmin}\{a\f:||a||_1=1,a\in\Lattice^\circ_+\},$
we again compute the set $\CP_\f(\succ,S) = \emph{\{(2,-1,1,-5,0,3)}\}$.
Note that $a\f=\frac{1}{4}(0,1,1,0,2,0)\cdot(0,1,0,5,0,0)=\frac{1}{4}$.
Then, for example,
$a\f'=\frac{1}{4}(0,1,1,0,2,0)\cdot\text{\emph{(1,-1,-1,-3,1,2)}}^+=
\frac{1}{2}>a\f$;
hence, the vector \emph{(1,-1,-1,-3,1,2)} is not needed in a $\f$-truncated
$\prec_c$-Gr\"obner basis of $\Lattice$.

On the other hand, if we used the vector
$a := \frac{1}{14}(0,0,1,5,0,8)\in\Lattice^\circ_+$,
then the quick truncation check $a\f < a\f'$ is useless, and we would compute
all five vectors of the $\prec_c$-Gr\"obner basis of $\Lattice$.

\item
Consider $\f:=(2,2,4,2,0,0)$ (see Figure \ref{Figure: Disconnected
Graph}a); $\FEASIBLE{\Lattice}{\f}=\{(2,1,1,1,2,1),(2,2,4,2,0,0)\}$.
Since, the feasible set consists of only two feasible solutions, the set
\emph{\{(0,1,3,1,-2,-1)\}} is a $\f$-truncated $\prec_c$-Gr\"obner basis of
$\Lattice$.

Using $\FEASIBLEX{\R}{\Lattice}{\f-\f'} \ne \emptyset$ as a check for
truncation, we compute the following $\f$-truncated $\prec_c$-Gr\"obner basis of
$\Lattice$:
$G := \text{\emph{\{(0,1,3,1,-2,-1),(1,0,2,-2,-1,1)\}}}.$
So, we have computed an additional vector \emph{(1,0,2,-2,-1,1)} that is not
strictly needed since $\FEASIBLE{\Lattice}{\f-\f'} = \emptyset$ even though
$\FEASIBLEX{\R}{\Lattice}{\f-\f'} \ne \emptyset$
where $\f':= \text{\emph{(1,0,2,-2,-1,1)}}^+=\text{\emph{(1,0,2,0,0,1)}}$.

Using the vector 
$a := \frac{1}{4}(0,1,0,1,0,2)=
\text{argmin}\{a\f:||a||_1=1,a\in\Lattice^\circ_+\}$
as a quick check for truncation,
we obtain the following $\f$-truncated $\prec_c$-Gr\"obner basis of
$\Lattice$:
\[G := \text{\emph{\{(1,-1,-1,-3,1,2), (1,0,2,-2,-1,1), (1,1,5,-1,-3,0),
(0,1,3,1,-2,-1)\}}}.\]
So here, we have computed three additional unnecessary vectors.
If instead we use the vector $\frac{1}{14}(0,0,1,5,0,8)$,
then we would have computed only three vectors.
In this case, there is no single vector in $a\in\Lattice^\circ_+$ that
results in only two vectors being computed.
\end{enumerate}
\end{example}

The following example demonstrates the potential speed increase from
computing a truncated Gr\"obner basis as opposed to computing a full Gr\"obner
basis.
We use three different methods for checking whether
$\f' \in \B{\Lattice}{\f}$ in order of increasing effectiveness:
\begin{enumerate}
\item $a\f < a\f'$ where
$a = \text{argmin} \{ a\f : ||a||_1 = 1, a \in \Lattice^\circ_+\}$,
\item $\FEASIBLEX{\R}{\Lattice}{\f-\f'} = \emptyset$, or
\item $\FEASIBLE{\Lattice}{\f-\f'} = \emptyset$.
\end{enumerate}
Since criterion (i) is in general much faster to check than (ii), we
always check (i) before (ii), and similarly, since criterion
(i) and (ii) are in general much faster to check than (iii), we always check
them both before applying criterion (iii).
We solve $\FEASIBLEX{\R}{\Lattice}{\f-\f'} = \emptyset$ using the simplex
algorithm implementation in the GLPK (GNU Linear Programming Kit)
package.\footnote{GLPK is open source and freely available from
\texttt{http://www.gnu.org/software/glpk/}}
We solve $\FEASIBLE{\Lattice}{\f-\f'} = \emptyset$ using the
branch-and-bound implementation in GLPK, which is not useful in practice,
but we use it here to show the sizes of minimal truncated Markov bases.

\begin{example}
\label{Example: GB}
Let $\Lattice = \Lattice_A := \{u : Au = \zero, u \in \Z^n\}$ where
\[
A =
\left[
\begin{array}{rrrrrrrrrrrrr}
15&  4& 14& 19&  2&  1& 10& 17& 11&  9&  4& 15& 20\\
18& 11& 13&  5& 16& 16&  8& 19& 18& 21&  5&  7&  1\\
11&  7&  8& 19& 15& 18& 14&  6&  1& 23& 11&  3& 10\\
17& 10& 13& 17& 16& 14& 15& 18&  3&  2&  1& 17&  1
\end{array}
\right]
\]
The size of a minimal Markov basis of $\Lattice$ is $10868$.
Let $c=(3,15,1,5,2,17,16,16,15,9,7,11,13)$. The size of a
minimal $\prec_c$-Gr\"obner basis of $\Lattice$ for some term order
$\prec$ is $24941$.
This takes $106.96$ seconds to compute using \emph{4ti2}.
In the following table, we list the time taken to compute truncated Gr\"obner
bases from the minimal Markov basis of $\Lattice$.
The first column lists the values used for $\f$.
In the following columns, we list the size of the computed set and the time
taken for each of the three possible ways to check whether
$\f'\in\B{\Lattice}{\f}$.
\begin{center}
\begin{tabular}{|c||c|c||c|c||c|c|}
\hline
$\f$ &
\multicolumn{2}{|c||}{$a\in\Lattice^\circ_+$} &
\multicolumn{2}{|c||}{$\FEASIBLEX{\R}{\Lattice}{\f-\f''}$} &
\multicolumn{2}{|c|}{$\FEASIBLE{\Lattice}{\f-\f''}$} \\
\hline
\hline
$(1,1,1,0,1,0,1,1,0,1,0,1,0)$&$307$&$0.02$s&$1$&$0.09$s&$0$&$0.10$s \\
\hline
$(1,0,1,0,3,0,1,5,0,1,0,9,0)$&$418$&$0.10$s&$36$&$0.19$s&$0$&$0.31$s\\
\hline
$(1,1,1,1,1,1,1,1,1,1,1,1,1)$&$6494$&$1.61$s&$201$&$1.20$s&$0$&$2.18$s \\
\hline
$(1,2,0,3,5,0,1,3,0,4,0,1,0)$&$12191$&$6.31$s&$5028$&$4.43$s&$158$&$186.81$s \\
\hline
$(19,7,3,8,13,11,1,15,4,8,17,9,5)$&$24748$&$108.19$s&$24334$&$107.25$s&$24284$&
$>3600$s\\
\hline
\end{tabular}
\end{center}
\end{example}

For the example above, choosing a single $a \in \Lattice^\circ_+$
where $a := \text{argmin} \{ a\f : ||a||_1 = 1, a \in \Lattice^\circ_+\}$ works
reasonably well when used for checking for truncation.
However, in general, using more than one such $a$ may be significantly better particularly when the support of $a$ (the set of non-zero components) is small.

Observe that to compute a truncated Gr\"obner bases in the previous example, we
first needed to compute a Markov basis, and in some cases, computing the
Markov basis took significantly longer than computing the truncated Gr\"obner
basis. This provides motivation for the next section in which we compute
truncated Markov basis.


%

\section{Computing truncated Markov bases}
\label{Section: Truncated Markov bases}
In this section, we give a Project-and-Lift algorithm for computing truncated
Markov bases.

Given $\sigma \subseteq \{1,\dots,n\}$, we define the projective map
$\pi_\sigma: \Z^n \mapsto \Z^{|\sigmabar|}$ that projects a vector in $\Z^n$
onto the $\sigmabar$ components.
We define $\Lattice^\sigma$ where $\sigma\subseteq\{1,\ldots,n\}$
as the projection of $\Lattice$ onto the $\sigmabar$
components -- that is, $\Lattice^\sigma = \pi_\sigma(\Lattice)$.
Note that $\Lattice^\sigma$ is also a lattice.
For ease of notation, we will often denote the singleton set $\{i\}$ as just
$\seti$, and so, for example, $\Lattice^{\{i\}}$ is denoted $\Lattice^\seti$ and
$\pi_{\{i\}}$ is denoted $\pi_\seti$.
It should be clear from the context
whether by $i$ we mean $\{i\}$ or just $i$.

The fundamental idea behind the Project-and-Lift algorithm is that,
for some $\f \in \Z^n$, using a set $S \subseteq \Lattice^\seti$ that is a
$\pi_\seti(\f)$-truncated Markov basis of $\Lattice^\seti$ for some
$i \in \{1,\dots,n\}$, we can compute a set $S'\subseteq \Lattice^\seti$ such
that $S'$ lifts to a $\f$-truncated Markov basis of $\Lattice$.
So, for some $\sigma \subseteq \{1,\dots,n\}$, since $\Lattice^\sigma$ is also a
lattice, starting with a $\f_\sigmabar$-truncated Markov basis of
$\Lattice^\sigma$, we can compute a $\pi_{\sigma\backslash\seti}(\f)$-truncated
Markov basis of $\Lattice^{\sigma\backslash\seti}$ for some $i \in \sigma$.
By doing this repeatedly for every $i \in \sigma$, we attain a $\f$-truncated
Markov basis of $\Lattice$.

First, we extend the definition of reduction paths and
$\f$-truncated Gr\"obner bases.
Given some vector $c \in \Q^n$, a path $(x^0,\ldots,x^k)$ in
$\GRAPH{\Lattice}{\f}{G}$ is an $c$-reduction path
if for all $j\in\{1,\ldots,k-1\}$,
we have either $c x^0 \ge c x^j$ or $c x^k \ge c x^j$.
A set $G \subseteq \Lattice$ is a
\textbf{$c$-Gr\"obner basis} of $\FEASIBLE{\Lattice}{\f}$ if
for every pair $x,y \in \FEASIBLE{\Lattice}{\f}$, there exists a
$c$-reduction path from $x$ to $y$ in $\GRAPH{\Lattice}{\f}{G}$.
A set $G \subseteq \Lattice$ is a $\f$-truncated
\textbf{$c$-Gr\"obner basis} of $\Lattice$ if for all $\f'\in\B{\Lattice}{\f}$,
$G$ is a $c$-Gr\"obner basis of $\FEASIBLE{\Lattice}{\f'}$.

The following lemma is fundamental to the Project-and-Lift algorithm.
Note that the property that $\ker(\pi_\seti)\cap\Lattice = \{\zero\}$
for some $i \in \{1,\dots,n\}$ means that the map $\pi_\seti$ from
$\Lattice$ to $\Lattice^\seti$ is a bijection and thus,
the inverse map $\pi^\inv_\seti:\Lattice^\seti \mapsto \Lattice$ is
well-defined (each vector in $\Lattice^\seti$ lifts to a unique vector in
$\Lattice$).
Moreover, by linear algebra, for all $u \in \Lattice$,
there must exist a vector $\omega^i \in \Q^{n-1}$ such that
$\omega^i \cdot \pi_\seti(u) = u_i$.
We always write such a vector as $\omega^i$.
Importantly, note that given $\f \in \Z^n$, we have
$\pi_\seti(\B{\Lattice}{\f}) \subseteq \B{\Lattice^\seti}{\pi_\seti(\f)}$
since if $\f' \in \B{\Lattice}{\f}$,
then there exists $\gamma \in \FEASIBLE{\Lattice}{\f-\f'}$, and so,
$\pi_\seti(\gamma)\in\FEASIBLE{\Lattice^\seti}{\pi_\seti(\f-\f')}\ne\emptyset$,
and consequently, $\pi_\seti(\f') \in \B{\Lattice^\seti}{\pi_\seti(\f)}$.

\begin{lemma}\label{Lemma:project-and-lift}
Let $i \in \{1,\dots,n\}$ where $\ker(\pi_\seti)\cap\Lattice = \{\zero\}$,
and let $S \subseteq \Lattice^\seti$. Let $\f \in \Z^n$.
If $S$ is a $\pi_\seti(\f)$-truncated $(\omegabar^i)$-Gr\"obner basis of
$\Lattice^\seti$, then $\pi^\inv_\seti(S)$ is a
$\f$-truncated $(\ebar^i)$-Gr\"obner basis of $\Lattice$.
\end{lemma}
\begin{proof}
Assume $S$ is a $\pi_\seti(\f)$-truncated $(\omegabar^i)$-Gr\"obner basis of
$\Lattice^\seti$.
Let $x,y \in \FEASIBLE{\Lattice}{\f'}$ for some $\f' \in \B{\Lattice}{\f}$.
We need to show that there is an $(\ebar^i)$-reduction path from $x$ to $y$ in
$\GRAPH{\Lattice}{\f'}{\pi^\inv_\seti(S)}$.
Let $\tilde{x}=\pi_\seti(x)$, $\tilde{y}=\pi_\seti(y)$, and
$\tilde{\f}'=\pi_\seti(\f') \in \B{\Lattice^\seti}{\pi_\seti(\f)}$.
By assumption, there exists an $(\omegabar^i)$-reduction path
$(\tilde{x}=\tilde{x}^0,\dots,\tilde{x}^k=\tilde{y})$
in $\GRAPH{\Lattice^\seti}{\tilde{\f}'}{S}$.
So, we have either $\omega^i\tilde{x}^j \ge \omega^i\tilde{x}$ or
$\omega^i\tilde{x}^j \ge \omega^i\tilde{y}$ for all $j$.
We now lift this $(\omegabar^i)$-reduction path in
$\GRAPH{\Lattice^\seti}{\tilde{\f}'}{S}$ to an $(\ebar^i)$-reduction
path in $\GRAPH{\Lattice}{\f'}{\pi^\inv_\seti(S)}$.
Let $x^j=\f'+\pi^\inv_\seti(\tilde{x}^j-\tilde{\f}')$
for all $j=0,\dots,k$.
Hence, $\pi_\seti(x^j)=\tilde{x}^j$ and $\omega^i\tilde{x}^j = x^j_i$,
and therefore, either $x^j_i \ge x_i$ or $x^j_i \ge y_i$.
Also, $x^j-x^{j-1} = \pi^\inv_\seti(\tilde{x}^j-\tilde{x}^{j-1}) \in
\pi^\inv_\seti(S)$ for all $j=1,\dots,k$.
Therefore, $(x=x^0,\dots,y^k=y)$ is an
$(\ebar^i)$-reduction path in
$\GRAPH{\Lattice}{\f}{\pi^\inv_\seti(S)}$ as required.
\end{proof}
The converse of Lemma \ref{Lemma:project-and-lift} is not true:
if $\pi^\inv_\seti(S)$ is a
$\f$-truncated $(\ebar^i)$-Gr\"obner basis of $\Lattice$,
then $S$ is \emph{not necessarily} a $\pi_\seti(\f)$-truncated 
$(\omegabar^i)$-Gr\"obner basis of $\Lattice$.
If $\pi_\sigma(\B{\Lattice}{\f}) = \B{\Lattice^\sigma}{\f_\sigmabar}$,
then the converse holds, since
if $\pi^\inv_\seti(S)$ is a
$(\ebar^i)$-Gr\"obner basis of $\FEASIBLE{\Lattice}{\f'}$ for some 
$\f' \in \Z^n$, then $S$ is a $(\omegabar^i)$-Gr\"obner basis of
$\FEASIBLE{\Lattice}{\bar{\f}'}$ where $\bar{\f}'=\pi_\seti(\f')$
because $(\ebar^i)$-reduction paths in
$\GRAPH{\Lattice}{\f'}{\pi^\inv_\seti(S)}$
project to $(\omegabar^i)$-reduction paths in
$\GRAPH{\Lattice}{\bar{\f}'}{S}$.
In general however,
$\pi_\seti(\B{\Lattice}{\f}) \subseteq \B{\Lattice^\seti}{\pi_\seti(\f)}$,
and we may have
$\pi_\seti(\B{\Lattice}{\f}) \subsetneq \B{\Lattice^\seti}{\pi_\seti(\f)}$.

By definition, an $(\ebar^i)$-Gr\"obner basis of $\Lattice$ is a Markov basis
of $\Lattice$. Conversely, a Markov basis of $\Lattice$ is also an
$(\ebar^i)$-Gr\"obner basis of $\Lattice$.
This follows since, given a Markov basis of $\Lattice$,
for any $x,y \in \FEASIBLE{\Lattice}{\f}$ for any $\f$,
there must exist a path from $x-\gamma$ to $y-\gamma$
where $\gamma = \min\{x_i,y_i\}\cdot e^i$, and by translating such a path by
$\gamma$, we get an $(\ebar^i)$-reduction path from $x$ to $y$.
So, we arrive at the following corollary.

\begin{corollary}
\label{Corollary: lifting a Markov basis}
Let $i \in \{1,\dots,n\}$ where $\ker(\pi_\seti)\cap\Lattice = \{\zero\}$,
and let $S \subseteq \Lattice^\seti$.
If $S$ is a $\pi_\seti(\f)$-truncated $(\omegabar^i)$-Gr\"obner basis of
$\Lattice^\seti$, then $\pi^\inv_\seti(S)$ is a $\f$-truncated 
Markov basis of $\Lattice$.
\end{corollary}

Given any vector $c \in \Z^n$ and a term order $\prec$ for $\Lattice$,
recall that for the ordering $\prec_c$, we have $x \prec_c y$ if $c x < c y$
or $c x = c y$ and $x \prec y$.
Also, recall that the order $\prec_c$ is a term order if and only if
$\IP{\Lattice}{c}{\zero} = 0$.
Importantly then, a $\prec_c$-reduction path is also a 
$c$-reduction path. So, we can compute a $\f$-truncated
$c$-Gr\"obner basis by computing a $\f$-truncated $\prec_c$-Gr\"obner
basis.

If $\IP{\Lattice}{\ebar^i}{\zero}:= 
\max\{x_i: x \in \Lattice, x \in \N^n\}=0$, 
we say that $i$ is \textbf{bounded} for $\Lattice$ and \textbf{unbounded}
otherwise.
Thus, $\prec_{\ebar^i}$ is a term order for $\Lattice$ if and only if $i$ is
bounded, and moreover, $\prec_{\omegabar^i}$ is a term order for
$\Lattice^\seti$ if and only if $i$ is bounded
since
\[\IP{\Lattice}{\ebar^i}{\zero} =
\max\{x_i: x \in \Lattice, x \in \N^n\} =
\max\{\omega^i x: x \in \Lattice^i, x \in \N^{n-1}\} =
\IP{\Lattice^\seti}{\omegabar^i}{\zero}.\]

Now if $i$ is bounded, then the ordering $\prec_{\omegabar^i}$ is a term order
for $\Lattice^\seti$, and so, given a set $S \subseteq \Lattice^\seti$ that is a
$\f$-truncated Markov basis of $\Lattice^\seti$, we can compute a
$\f$-truncated $(\omegabar^i)$-Gr\"obner basis of $\Lattice^\seti$
using Algorithm \ref{Algorithm: Completion procedure}. In other words,
the set $S' = \CP_{\pi_\seti(\f)}(\prec_{\omegabar^i},S)$ is a
$\pi_\seti(\f)$-truncated $(\omegabar^i)$-Gr\"obner basis of $\Lattice^\seti$,
and by Corollary \ref{Corollary: lifting a Markov basis},
the set $\pi^\inv_\seti(S')$ is a $\f$-truncated Markov basis of $\Lattice$.

If $i$ is unbounded, then computing a $\f$-truncated Markov basis of
$\Lattice$ from a $\pi_\seti(\f)$-truncated Markov basis of $\Lattice^\seti$
is actually more straight-forward than otherwise.
Crucially, $i$ is unbounded if and only if $\IP{\Lattice}{\ebar^i}{\zero}>0$ or
equivalently there exists $u \in \Lattice \cap \N^n$ where $u_i > 0$. 
Then assuming $\ker(\pi_\seti)\cap \Lattice = \emptyset$,
given a set $S \subseteq \Lattice^\seti$ that is a 
$\pi_\seti(\f)$-truncated Markov basis of $\Lattice^\seti$,
it suffices to add $u$ to $\pi^\inv_\seti(S')$ to create a $\f$-truncated
Markov basis of $\Lattice$
(see Lemma \ref{Lemma: gen non-negative vector} below).
We can use linear programming to check whether $i$ is unbounded and also to
find such a $u \in \Lattice\cap\N^n$ where $u_i > 0$.
\begin{lemma}
\label{Lemma: gen non-negative vector}
Let $i \in \{1,\dots,n\}$ and $\ker(\pi_\seti)\cap\Lattice = \emptyset$,
and $u \in \Lattice \cap \N^n$ where $u_i > 0$.
If $S \subseteq \Lattice^\seti$ is a $\pi_\seti(\f)$-truncated Markov basis
of $\Lattice^\seti$, then $\pi^\inv_\seti(S)\cup\{u\}$ is a
$\f$-truncated Markov basis of $\Lattice$.
\end{lemma}
\begin{proof}
Let $x',y' \in \FEASIBLE{\Lattice}{\f'}$ for some $\f' \in \B{\Lattice}{\f}$.
Since $S$ is a $\pi_\seti(\f)$-truncated Markov basis of $\Lattice^\seti$,
there exists a path from $\pi_\seti(x')$ to $\pi_\seti(y')$.
We can convert this path into a $(\omegabar^i)$-reduction path
in $\GRAPH{\Lattice^\seti}{\pi_\seti(\f')}{S\cup\pi_\seti(u)}$
by adding $\pi_\seti(u)$ to the start of the path as many times as
necessary and subtracting $\pi_\seti(u)$ from the end of the path the same
number of times.
This works since $\omega^i(\pi_\seti(u))=u_i > 0$ and $u \ge 0$.
As in Lemma \ref{Lemma:project-and-lift}, 
we can then lift this to a path from $x'$ to $y'$ in 
$\GRAPH{\Lattice}{\f'}{\pi^\inv_\seti(S)\cup\{u\}}$.
\end{proof}

We can apply the above reasoning to compute a Markov basis of
$\Lattice^{\sigma\setminus\seti}$ from a Markov basis of $\Lattice^\sigma$ for
some $\sigma \subseteq \{1,\dots,n\}$ and $i \in \sigma$.
First, analogously to $\pi_\seti$ and $\prec_{\omegabar^i}$ in the context
of $\Lattice^\seti$ and $\Lattice$, we define
$\pi^\sigma_\seti$ and $\prec^\sigma_{\omegabar^i}$ in the same way except
in the context of $\Lattice^\sigma$ and $\Lattice^{\sigma\setminus\seti}$
respectively.

We can now present our Project-and-Lift algorithm
(Algorithm \ref{Project-and-Lift algorithm}).

\begin{algorithm} 
\caption{Project-and-Lift algorithm}
\begin{algorithmic}
\REQUIRE a lattice $\Lattice$ and a vector $\f \in \Z^n$.
\ENSURE a $\f$-truncated Markov basis $M$ of $\Lattice$
\STATE Find a set $\sigma\subseteq\{1,\ldots,n\}$ such that
$\ker(\pi_\sigma) \cap \Lattice = \{\zero\}$.
\STATE Compute a set $M \subseteq \Lattice^\sigma$ that is a
$\f_\sigmabar$-truncated Markov basis of $\Lattice^\sigma$.
\WHILE{$\sigma\neq\emptyset$}
  \STATE Select $i \in \sigma$
  \IF{$i$ is bounded}
     \STATE $G:= \CP_{\f_\sigmabar}(\prec^\sigma_{\omegabar^i},M)$
     \STATE $M:= (\pi^\sigma_\seti)^\inv(G)$
  \ELSE
     \STATE Compute $u \in \Lattice^{\sigma\setminus\seti} \cap \N^n$
            such that $u_i > 0$
     \STATE $M:= (\pi^\sigma_\seti)^\inv(M) \cup \{u\}$
  \ENDIF
  \STATE $\sigma:=\sigma\setminus\seti$
  \STATE $M:= \{ u \in M: u^+ \in \B{\Lattice^\sigma}{\f_\sigmabar}\}$
\ENDWHILE
\STATE \Kw{return} $M$.
\end{algorithmic}
\label{Project-and-Lift algorithm}
\end{algorithm}
\begin{lemma}
\label{Lemma: Project-and-Lift algorithm terminates and is correct}
Algorithm \ref{Project-and-Lift algorithm} terminates and satisfies its
specifications.
\end{lemma}
\begin{proof}
Algorithm \ref{Project-and-Lift algorithm} terminates, since
Algorithm \ref{Algorithm: Completion procedure}, which computes
$\CP_{\f_\sigmabar}(\prec^\sigma_{\omegabar^i},M))$, always terminates.

We claim that for each iteration of the algorithm, $M$ is a
$\f_\sigmabar$-truncated Markov basis of
$\Lattice^\sigma$ and $\ker(\pi_\sigma) \cap \Lattice = \{\zero\}$;
therefore, at termination, $M$ is a $\f$-truncated Markov basis of $\Lattice$.
This is true for the first iteration, so we assume it is true for
the current iteration.

If $\sigma = \emptyset$, then there is nothing left to do, so assume otherwise.
Since by assumption,
$\ker(\pi_\sigma) \cap \Lattice^\sigma = \{\zero\}$, we must have
$\ker(\pi^\sigma_\seti) \cap \Lattice^{\sigma\setminus\seti} = \{\zero\}$,
and so, the inverse map
$(\pi^\sigma_\seti)^\inv: \Lattice^\sigma \rightarrow
\Lattice^{\sigma\setminus\seti}$ is well-defined.
Let $i \in \sigma$, and $\sigma' := \sigma\setminus\seti$.
If $i$ is bounded, then let
$G:=\CP_{\f_\sigmabar}(\prec^\sigma_{\omegabar^i},M)$;
then, $G$ is a $\prec^\sigma_{\omegabar^i}$-Gr\"obner basis of
$\Lattice^\sigma$.
Let $M':=(\pi^\sigma_\seti)^\inv(G)$,
and then by Corollary \ref{Corollary: lifting a Markov basis},
$M'$ is a $\pi_{\sigma'}(\f)$-truncated Markov basis of $\Lattice^{\sigma'}$.
Otherwise, let $M:= (\pi^\sigma_\seti)^\inv(M) \cup \{u\}$ where
$u \in \Lattice^{\sigma\setminus\seti} \cap \N^n$
such that $u_i > 0$, and by Lemma \ref{Lemma: gen non-negative vector},
$M'$ is a $\pi_{\sigma'}(\f)$-truncated Markov basis of $\Lattice^{\sigma'}$.
Also, $M'' := \{ u \in M': u^+ \in \B{\Lattice^{\sigma'}}{\pi_{\sigma'}(\f)}\}$
must be a $\pi_{\sigma'}(\f)$-truncated Markov basis of $\Lattice^{\sigma'}$.

Lastly, since $\sigma' \subseteq \sigma$,
we must have $\ker(\pi_{\sigma'}) \cap \Lattice = \{\zero\}$.
Thus, the claim is true for the next iteration.
\end{proof}
Initially in our Project-and-Lift algorithm, we need to find a set
$\sigma\subseteq\{1,\ldots,n\}$
such that $\ker(\pi_\sigma) \cap \Lattice^\sigma = \{\zero\}$,
and then, we need to compute a Markov basis for $\Lattice^\sigma$.
This is actually quite straight-forward and can be done in polynomial time.
Let $B$ be a basis for the lattice $\Lattice$ ($\Lattice$ is spanned by the rows
of the matrix $B$).
Let $k := \rank(B)$.
Any $k$ linearly independent columns of $B$ then suffice to
give a set $\sigmabar$ such that every vector in $\Lattice^\sigma$
lifts to a unique vector in $\Lattice$; that is,
$\ker(\pi_\sigma) \cap \Lattice^\sigma = \{\zero\}$.
Such a set $\sigma$ can be found via Gaussian elimination.
Let $S = \pi_\sigma(B)$; then, $S$ spans $\Lattice^\sigma$, and
$S \in \Z^{k\times k}$ since $|\sigmabar|=k$.
Let $S'$ be an upper triangle matrix with positive diagonal entries and
non-positive entries elsewhere such
that (the rows of) $S'$ span $\Lattice^\sigma$.
We can always construct such a matrix $S'$ from $S$
in polynomial time using the \emph{Hermite Normal Form} (HNF) algorithm 
(see for example \cite{NemhauserWolsey1988}).
Also, $S'$ is a Markov basis of $\Lattice^\sigma$ since it is actually a
Gr\"obner basis of $\Lattice^\sigma$ with respect to a \emph{lexicographic}
ordering, and thus, a Markov basis of $\Lattice^\sigma$.

\begin{example}
Consider again the set $S := $\emph{\{(1,-1,-1,-3,1,2),(1,0,2,-2,-1,1)\}}, and
the lattice $\Lattice \subseteq \Z^6$ spanned by $S$.
Let again $\f:=(0,1,0,5,0,0)$
(see Figure \ref{Figure: Disconnected Graph}b).
Recall that the minimal $\f$-truncated $\prec_c$-Gr\"obner basis
is \emph{\{(2,-1,1,-5,0,3)}\}, since
$\FEASIBLE{\Lattice}{\f}=\{(2,0,1,0,0,3),(0,1,0,5,0,0)\}$.

Let $\sigma =\{3,4,5,6\}$.
Then, $\ker(\pi_\sigma) \cap \Lattice^\sigma = \{\zero\}$. 
Note that $\pi_\sigma(S)=$\emph{\{(1,-1),(1,0)\}}.
The set $M = \{(1,0),(0,1)\}$ is a Markov basis of $\Lattice^\sigma$.
\begin{enumerate}
\item Set $i := 3$. Then, $i$ is unbounded since
$(1,0,2)\in\Lattice^{\sigma\setminus\seti}$.
Thus, $M = \{(1,0,2),(0,1,3)\}$ is a Markov basis of
$\Lattice^{\sigma\setminus\seti}$.
Set $\sigma = \{4,5,6\}$.
\item Set $i := 5$. Then, $i$ is bounded.
$M = $\emph{\{(-2,1,-1),(1,0,2)\}} is a minimal $(0,1,0)$-truncated
$\prec_{\omegabar^i}$-Gr\"obner basis of $\Lattice^\sigma$.
Note that $\f_\sigmabar = (0,1,0)$.
So, $M = $\emph{\{(-2,1,-1,0),(1,0,2,-1)\}} is a $(0,1,0,0)$-truncated
Markov basis of $\Lattice^{\sigma\setminus\seti}$.
Set $\sigma = \{4,6\}$.
\item Set $i := 4$. Then, $i$ is bounded.
$M = $\emph{\{(2,-1,1,0)\}} is a minimal $(0,1,0,0)$-truncated
$\prec_{\omegabar^i}$-Gr\"obner basis of $\Lattice^\sigma$.
So, $M = $\emph{\{(2,-1,1,-5,0)\}} is a $(0,1,0,5,0)$-truncated Markov basis
of $\Lattice^{\sigma\setminus\seti}$.
Set $\sigma = \{6\}$.
\item Set $i := 6$. Then, $i$ is bounded.
$M = $\emph{\{(-2,1,-1,5,0)\}} is a minimal $(0,1,0,5,0)$-truncated
$\prec_{\omegabar^i}$-Gr\"obner basis of $\Lattice^\sigma$.
So, $M = $\emph{\{(-2,1,-1,5,0,-3)\}} is a $(0,1,0,5,0,0)$-truncated 
Markov basis of $\Lattice$.
\end{enumerate}
In the above Markov basis computation, the size of the set $M$ was never
larger than $2$ although the size of a full minimal Markov basis is $5$.
\end{example}

In the next example, we show the computational benefits of computing a
truncated Markov basis as opposed to the full Markov basis.

\begin{example}
Let $\Lattice = \Lattice_A$ for the matrix $A$ given in Example
\ref{Example: GB}.
The size of a minimal Markov basis of $\Lattice$ is $10868$.
It takes $36.39$ seconds to compute.

Let $\f = (1,1,1,1,1,1,1,1,1,1,1,1,1)$.  The size of a minimal $\f$-truncated
Markov basis of $\Lattice$ is 0. In the next table, we list the times taken to
compute a $\f$-truncated Markov basis using the three different criteria for truncation.
We also list the sizes of the intermediate Gr\"obner basis computations.
The number at the top of the column is the size of $\sigma$ as used in algorithm
\ref{Project-and-Lift algorithm}.
The first row is the case without using truncation.
Note how the intermediate sizes of the truncated computations remain much
smaller than the final size of a full Markov basis of $\Lattice$.
\begin{center}
\begin{tabular}{|c|c|c|c|c|c|c|}
\hline
Truncation &$4$&$3$&$2$&$1$&$0$&Time\\
\hline
\hline
none  & $545$ & $1822$ &$3681$&$12573$&$10868$& $36.39$s \\
\hline
$a\in\Lattice^\circ_+$  & $545$ & $977$ &$1302$&$1846$&$564$& $1.20$s \\
\hline
$\FEASIBLEX{\R}{\Lattice}{\f-\f''}$ & $545$ & $977$ &$878$ &$697$ &$194$& $2.16$s \\
\hline
$\FEASIBLE{\Lattice}{\f-\f''}$    & $545$ &   $6$ &  $3$ &  $4$ & $0$ &$20.73$s \\
\hline
\end{tabular}
\end{center}

In the following table, we list the times taken to compute a minimal truncated
Markov basis for different $\f$.
The first column lists the values used for $\f$.
In the next columns, we list the size of the computed set and the time taken for
each of the three possible ways to check whether $\f' \in \B{\Lattice}{\f}$.
\begin{center}
\begin{tabular}{|c||c|c||c|c||c|c|}
\hline
$\f$ &
\multicolumn{2}{|c||}{$a\in\Lattice^\circ_+$} &
\multicolumn{2}{|c||}{$\FEASIBLEX{\R}{\Lattice}{\f-\f''}$} &
\multicolumn{2}{|c|}{$\FEASIBLE{\Lattice}{\f-\f''}$} \\
\hline
\hline
$(1,1,1,0,1,0,1,1,0,1,0,1,0)$&   $4$& $0.76$s&   $1$& $0.90$s&  $0$&$3.33$s \\
\hline
$(1,0,1,0,3,0,1,5,0,1,0,9,0)$& $158$&$16.70$s&  $36$&$21.28$s&  $0$&$>3600$s\\
\hline
$(1,1,1,1,1,1,1,1,1,1,1,1,1)$& $546$& $1.20$s& $194$& $2.16$s&  $0$&$20.73$s \\
\hline
$(1,2,0,3,5,0,1,3,0,4,0,1,0)$&$7381$& $7.61$s&$3734$&$11.47$s&$146$&$>3600$s \\
\hline
$(19,7,3,8,13,11,1,15,4,8,17,9,5)$&$10814$&$39.04$s&$10761$&$44.07$s&$10739$&$>3600$s \\
\hline
\end{tabular}
\end{center}
Not only have we reduced the time to compute a Markov basis, we also have
reduced the time to compute a truncated Gr\"obner basis by using a truncated 
Markov basis instead of a full Markov basis.
The following table lists the times for computing a truncated Gr\"obner basis
from a truncated Markov basis using the same cost vector
$c=(3,15,1,5,2,17,16,16,15,9,7,11,13)$ as before.
\begin{center}
\begin{tabular}{|c||c|c||c|c||c|c|}
\hline
$\f$ &
\multicolumn{2}{|c||}{$a\in\Lattice^\circ_+$} &
\multicolumn{2}{|c||}{$\FEASIBLEX{\R}{\Lattice}{\f-\f''}$} &
\multicolumn{2}{|c|}{$\FEASIBLE{\Lattice}{\f-\f''}$} \\
\hline
\hline
$(1,1,1,0,1,0,1,1,0,1,0,1,0)$&    $4$&$0.00$s&$1$&$0.00$s&$0$&$0.00$s \\
\hline
$(1,0,1,0,3,0,1,5,0,1,0,9,0)$&  $167$&$0.00$s&$36$&$0.19$s&$0$&$0.00$s\\
\hline
$(1,1,1,1,1,1,1,1,1,1,1,1,1)$&  $844$&$0.05$s&$201$&$0.11$s&$0$&$0.00$s \\
\hline
$(1,2,0,3,5,0,1,3,0,4,0,1,0)$&$11768$&$5.92$s&$5028$&$4.07$s&$158$&$5.98$s \\
\hline
$(19,7,3,8,13,11,1,15,4,8,17,9,5)$& $24729$ & $109.57$ & $24334$ & $107.18$s &
$24284$ & $>3600$s \\
\hline
\end{tabular}
\end{center}
\end{example}

Observe that in the previous example, the size of a truncated Markov basis for
the fiber where $\f = (1,0,1,0,3,0,1,5,0,1,0,9,0)$ was very small, but it took
much longer to compute than a much larger truncated Markov basis of other
fibers.
This anomaly can be attributed to the order in which the variables are chosen
during the project-and-lift algorithm. If we reorder the variables so that the
zero components in $\f = (1,0,1,0,3,0,1,5,0,1,0,9,0)$ are chosen first, then the
algorithm computes a truncated Markov basis much faster.
Hence, the efficiency of the algorithm is sensitive to the order in which the
variables are chosen; therefore, future experimentation is needed to determine a
heuristic for choosing a \emph{good} variable ordering.

\section{Feasibility} \label{Section:Feasibility}
The project-and-lift algorithm can be used to find a feasible solution of
$\FEASIBLE{\Lattice}{\f}$ for any $\f \in \Z^n$.
Alternatively, in the truncated case, given a particular $\f$, the
project-and-lift algorithm can be used to find a feasible solution of
$\FEASIBLE{\Lattice}{\f}$ and also a feasible solution of
$\FEASIBLE{\Lattice}{\f'}$ for any $\f' \in \B{\Lattice}{\f}$.
Thus, at the same time as finding a Markov basis, the algorithm can also compute
a feasible solution.


The basic idea is that, given some $\f \in \Z^n$ and a feasible solution of
$\FEASIBLE{\Lattice^\seti}{\pi_\seti(\f)}$ for some $i \in \{1,\dots,n\}$
where $\ker(\pi_\seti)\cap\Lattice = \{\zero\}$, we can
construct a feasible solution of $\FEASIBLE{\Lattice}{\f}$ if such a feasible
solution exists.
Hence, for some $\sigma \subseteq \{1,\dots,n\}$, starting with a feasible solution
of $\FEASIBLE{\Lattice^\sigma}{\f_\sigmabar}$ (we choose $\sigma$ such that
this feasible solution is easy to find), we can compute a feasible solution of
$\FEASIBLE{\Lattice^{\sigma\backslash\seti}}{\pi_{\sigma\backslash\seti}(\f)}$
for some $i \in \sigma$.  By doing this repeatedly for every $i \in \sigma$,
we attain a feasible solution of $\FEASIBLE{\Lattice}{\f}$.

We now show how to construct a feasible solution of $\FEASIBLE{\Lattice}{\f}$
from a feasible solution of $\FEASIBLE{\Lattice^\seti}{\pi_\seti(\f)}$ for some
$i \in \{1,...,n\}$ where $\ker(\pi_\seti)\cap\Lattice = \{\zero\}$.
Let $\f \in \Z^n$, and $x \in \FEASIBLE{\Lattice^\seti}{\pi_\seti(\f)}$.
First, we lift $x$.
Let $x' := \pi^\inv_\seti(x-\pi_\seti(\f)) + \f$.
Note that $(x-\pi_\seti(\f)) \in \Lattice^\seti$ and so $\pi^\inv_\seti$ is
well-defined.
Now, we have $\pi_\seti(x') = x$ and so $x'_j \ge 0$ for
all $j \ne i$.
Also, $x'\in\FEASIBLEX{\Z}{\Lattice}{\f}$.
If $x'_i \ge 0$, then $x' \in \FEASIBLE{\Lattice}{\f}$ and we are done.
So, assume $x'_i < 0$.
Now, if $i$ is unbounded, then there exists $u \in \Lattice \cap \N^n$ where
$u_i > 0$; therefore, $x'+\lambda u$ for some $\lambda \in \N$ is
non-negative on the $i$th component and thus a feasible solution
of $\FEASIBLE{\Lattice}{\f}$.
If $i$ is bounded, then we can compute a set $G$ that is a
$\pi_\seti(\f)$-truncated $\prec_{\omegabar^i}$-Gr\"obner basis of
$\Lattice^\seti$ for some term order $\prec$.
Let $x'' = \NF(x,G)$. Hence, $x''$ is the optimal solution of 
$\IP{\Lattice^\seti}{\prec_{\omegabar^i}}{\pi_\seti(\f)}$
and so also an optimal solution of
$\IP{\Lattice^\seti}{\omegabar^i}{\pi_\seti(\f)}$.
Thus conceptually, when computing $\NF(x,G)$,
we are just maximising the $i$th component which will thus become non-negative if a
feasible solution of $\FEASIBLE{\Lattice}{\f}$ exists.
Therefore, from Lemma \ref{Lemma:feasibility} below,
either $x' := \pi^\inv_\seti(x''-\pi_\seti(\f)) + \f$ is a feasible
solution of $\FEASIBLE{\Lattice}{\f}$ or $\FEASIBLE{\Lattice}{\f} = \emptyset$.
\begin{lemma}\label{Lemma:feasibility}
Let $i \in \{1,\dots,n\}$ where $i$ is bounded and
$\ker(\pi_\seti)\cap\Lattice = \{\zero\}$.
Let $\f \in \Z^n$, and let $x \in \FEASIBLE{\Lattice^\seti}{\pi_\seti(\f)}$ be a
$(\omegabar^i)$-minimal solution of $\FEASIBLE{\Lattice^\seti}{\pi_\seti(\f)}$.
Then, $\FEASIBLE{\Lattice}{\f} \ne \emptyset$ if and only if
$(\pi^\inv_\seti(x-\pi_\seti(\f))+\f) \in \FEASIBLE{\Lattice}{\f}$.
\end{lemma}
\begin{proof}
Let $x' = (\pi^\inv_\seti(x-\pi_\seti(\f))+\f)$.
If $x' \in \FEASIBLE{\Lattice}{\f}$, then
$\FEASIBLE{\Lattice}{\f} \ne \emptyset$ by definition. We now prove the
converse. Assume $\FEASIBLE{\Lattice}{\f} \ne \emptyset$, and let
$x'' \in \FEASIBLE{\Lattice}{\f}$.
Then, $\pi_\seti(x'') \in \FEASIBLE{\Lattice^\seti}{\pi_\seti(\f)}$ and so
$\omegabar^i\pi_\seti(x'') \ge \omegabar^i x$ which implies that
$x''_i \le (\pi^\inv_\seti(x''))_i$. Therefore, $\pi^\inv_\seti(x'')$ is
non-negative and a feasible solution of $\FEASIBLE{\Lattice}{\f}$.
\end{proof}

To compute a feasible solution of $\FEASIBLE{\Lattice}{\f}$, we need to start
from a feasible solution of $\FEASIBLE{\Lattice^\sigma}{\f_\sigmabar}$ for
some $\sigma \subseteq \{1,...,n\}$ where
$\ker(\pi_\sigma) \cap \Lattice = \{\zero\}$.
As before in Section \ref{Section: Truncated Markov bases}, we can find a set
$\sigma \subseteq \{1,..,n\}$ such that 
$\ker(\pi_\sigma) \cap \Lattice^\sigma = \{\zero\}$.
Moreover, we can also find a set $S \subseteq \Lattice^\sigma$ such that 
$S$ spans $\Lattice^\sigma$, and $S$ is an upper triangle square matrix with
positive diagonal entries and non-positive entries elsewhere.
Now, the vector $\f_\sigmabar$ is a solution to the relaxation
$\FEASIBLEX{\Z}{\Lattice^\sigma}{\f_\sigmabar}$ (non-negativity constraints
are removed).
Then, we can add appropriate non-negative multiples of the vectors in $S$ to
$\f_\sigmabar$ such that it becomes non-negative, and thus, we arrive at a
feasible solution of $\FEASIBLE{\Lattice^\sigma}{\f_\sigmabar}$.

\begin{algorithm}[ht]
\caption{Feasibility algorithm}
\begin{algorithmic}
\REQUIRE a lattice $\Lattice$ and a vector $\f \in \Z^n$.
\ENSURE a feasible solution $x \in \FEASIBLE{\Lattice}{\f}$ or infeasible.
\STATE Find a set $\sigma\subseteq\{1,\ldots,n\}$ such that
$\ker(\pi_\sigma) \cap \Lattice = \{\zero\}$.
\STATE Compute a set $M \subseteq \Lattice^\sigma$ that is a
$\f_\sigmabar$-truncated Markov basis of $\Lattice^\sigma$.
\STATE Compute a feasible solution $x\in\FEASIBLE{\Lattice^\sigma}{\f_\sigmabar}$.
\WHILE{$\sigma\neq\emptyset$}
  \STATE Select $i \in \sigma$
  \IF{$i$ is bounded}
     \STATE $G:= \CP_{\f_\sigmabar}(\prec^\sigma_{\omegabar^i},M)$
     \STATE $M:= (\pi^\sigma_\seti)^\inv(G)$
     \STATE $x:= \NF(x,G)$
     \STATE $x:= (\pi^\sigma_\seti)^\inv(x - \f_\sigmabar) +
     \pi_{\sigma\setminus\seti}(\f)$
     \STATE \Kw{if} $x_i < 0$ \Kw{then} \Kw{return} \emph{infeasible} 
  \ELSE
     \STATE Compute $u \in \Lattice^{\sigma\setminus\seti} \cap \N^n$
            such that $u_i > 0$
     \STATE $M:= (\pi^\sigma_\seti)^\inv(M) \cup \{u\}$
     \STATE $x:= (\pi^\sigma_\seti)^\inv(x - \f_\sigmabar) +
     \pi_{\sigma\setminus\seti}(\f)$
     \STATE $x:=x + \lambda u$ where $\lambda \in \N$ such that
       $(x + \lambda u)_i > 0$
  \ENDIF
  \STATE $\sigma:=\sigma\setminus\seti$
  \STATE $M:= \{ u \in M: u^+ \in \B{\Lattice^\sigma}{\f_\sigmabar}\}$
\ENDWHILE
\STATE \Kw{return} $x$.
\end{algorithmic}
\label{Algorithm: Feasibility}
\end{algorithm}

See Algorithm \ref{Algorithm: Feasibility} for a description of the
feasibility algorithm.

Finally, note that each individual step needed to compute a feasible solution of
$\FEASIBLE{\Lattice}{\f}$ is performed during the project-and-lift algorithm
\ref{Project-and-Lift algorithm}.
So, at the same time as computing a Markov basis, we can compute a feasible
solution.
Moreover, we can compute feasible solutions for many different fibers
simultaneously, and thus, we can avoid repeating the same Gr\"obner basis
computations for each fiber.

\begin{example}
We apply the above method to find a feasible solution of equality
constrained integer knapsack problems (see \cite{AardalLenstra2004}):
\[ \Feasible = \{x : A x = b, x \in \N^n\}\]
where $A \in \N^{1 \times n}$ and $b \in \N$.
Let $\Lattice = \Lattice_A := \{u:Au = \zero, u \in \Z^n\}$.
Then $\Feasible = \FEASIBLE{\Lattice}{\f}$ where
$\f \in \{x : A x = b, x \in \Z^n\}$.
Finding such a vector $\f$ can be done in polynomial time using the
HNF algorithm.
If no such $\f$ exists, then the original problem $\Feasible$ is infeasible.
Computing a Markov basis and thus solving the feasibility problem of any
such knapsack problem involves only one Gr\"obner basis
computation.\footnote{This Gr\"obner basis method for computing a feasible
solution for the special case of equality constrained integer knapsack problems
was found independently by Bjarke H. Roune (\cite{Roune06}).}

Consider the following knapsack feasibility problem:
\[\Feasible:=\{x:12223x_1+12224x_2+36674x_3+61119x_4+85569x_5=89643481:x\in\N^5\}.\]
Let $\Lattice = \Lattice_A$ where
$A= \begin{bmatrix} 12223 & 12224 & 36674 & 61119 & 85569 \end{bmatrix}$.
Then, the set
\[S =
\{\text{\emph{(12224,-12223,0,0,0),(2,-5,1,0,0),(-1,-4,0,1,0),(1,-8,0,0,1)}}\}\]
spans $\Lattice$.
Let $\f = (-4889,12222,0,0,0)$, then $\FEASIBLE{\Lattice}{\f} = \Feasible$.
Let $\sigma = \{1\}$. Then, $\ker(\pi_\sigma)\cap \Lattice = \{\zero\}$, and
also, we have
$\f_\sigmabar = x = (12222,0,0,0)\in\FEASIBLE{\Lattice^\sigma}{\f_\sigmabar}$.
A Markov basis of $\Lattice^\sigma$ is
\[ S' = \text{
\emph{\{(1,0,0,-1528),(0,1,0,-7640),(0,0,1,-6112),(0,0,0,12223)\}}
}.\]
Recall that $\omega^i \in \Q^{n-1}$ such that
$\omega^iu = (\pi_\seti^\inv(u))_i$.
A $\prec_{\omegabar^i}$-Gr\"obner basis of $\Lattice^\sigma$ is
\[ G = \text{
\emph{\{(3,2444,0,0),(-2,2445,0,0),(5,-1,0,0),(-4,0,1,0),(-3,-1,0,1)\}}
}.\]
The normal form of $x$ is $\NF(x,G) = x' = (2,2444,0,0)$.
We then lift this back into the original space,
$x'' = (x'-\f_\sigmabar)+\f = (-1,2,2444,0,0)$.
This is not a feasible solution, and the problem is therefore infeasible.

We list the times to solve the feasibility problem in Figure
\ref{figure:knapsacks}.
In each case, the problem was infeasible.
The right hand sides used for each knapsack is the Frobenius number, that is,
the largest infeasible right hand side.
\emph{
\begin{figure}[ht]
\begin{center}
\begin{tabular}{|l|l|c|c|}
\hline
Problem & Equality constraint ($A$) & RHS ($b$) & Time \\
\hline
\hline
cuww1 & \small{12223 12224 36674 61119 85569} & \small{89643481} & 0.00s \\
\hline
cuww2 & \small{12228 36679 36682 48908 61139 73365} & \small{89716838} & 0.00s \\
\hline
cuww3 & \small{12137 24269 36405 36407 48545 60683} & \small{58925134} & 0.00s \\
\hline
cuww4 & \small{13211 13212 39638 52844 66060 79268 92482} & \small{104723595} & 0.00s \\
\hline
cuww5 & \small{13429 26850 26855 40280 40281 53711 53714 67141} &
\small{45094583} & 0.00s \\
\hline
prob1  &\small{25067 49300 49717 62124 87608 88025 113673 119169} &
\small{3367335} & 0.00s \\
\hline
prob2  &\small{11948 23330 30635 44197 92754 123389 136951 140745} &
\small{14215206} & 0.00s \\
\hline
prob3  &\small{39559 61679 79625 99658 133404 137071 159757 173977} &
\small{58424799} & 0.02s \\
\hline
prob4  &\small{48709 55893 62177 65919 86271 87692 102881 109765} &
\small{60575665} & 0.01s \\
\hline
prob5  &\small{28637 48198 80330 91980 102221 135518 165564 176049} &
\small{62442884} & 0.01s \\
\hline
prob6  &\small{20601 40429 42407 45415 53725 61919 64470 69340 78539 95043} &
\small{22382774} & 0.30s \\
\hline
prob7  &\small{18902 26720 34538 34868 49201 49531 65167 66800 84069 137179} &
\small{27267751} & 0.00s \\
\hline
prob8  &\small{17035 45529 48317 48506 86120 100178 112464 115819 125128 129688}
& \small{21733990} & 0.01s \\
\hline
prob9  &\small{13719 20289 29067 60517 64354 65633 76969 102024 106036 199930} &
\small{13385099} & 0.01s \\
\hline
prob10 &\small{45276 70778 86911 92634 97839 125941 134269 141033 147279 153525}
& \small{106925261} & 0.05s \\
\hline
prob11 &\small{11615 27638 32124 48384 53542 56230 73104 73884 112951 130204} &
\small{577134} & 0.48s \\
\hline
prob12 &\small{14770 32480 75923 86053 85747 91772 101240 115403 137390 147371}
& \small{944183} & 0.32s \\
\hline
prob13 &\small{15167 28569 36170 55419 70945 74926 95821 109046 121581 137695} &
\small{6765260} & 0.78s \\
\hline
prob14 &\small{11828 14253 46209 52042 55987 72649 119704 129334 135589 138360}
& \small{80230} & 0.23s \\
\hline
prob15 &\small{13128 37469 39391 41928 53433 59283 81669 95339 110593 131989} &
\small{1663281} & 0.17s \\
\hline
prob16 &\small{35113 36869 46647 53560 81518 85287 102780 115459 146791 147097}
& \small{109710} & 0.75s \\
\hline
prob17 &\small{14054 22184 29952 64696 92752 97364 118723 119355 122370 140050}
& \small{752109} & 0.22s \\
\hline
prob18 &\small{20303 26239 33733 47223 55486 93776 119372 136158 136989 148851}
& \small{783879} & 0.51s \\
\hline
prob19 &\small{20212 30662 31420 49259 49701 62688 74254 77244 139477 142101} &
\small{677347} & 0.29s \\
\hline
prob20 &\small{32663 41286 44549 45674 95772 111887 117611 117763 141840 149740}
& \small{1037608} & 0.45s \\ \hline
\end{tabular}
\end{center}
\caption{Hard Knapsack Constraint Instances.}
\label{figure:knapsacks}
\end{figure}
}
\end{example}
In the paper \cite{AardalLenstra2004}, the feasibility is solved
problem for the same set of equality constrained integer knapsack problems by
using a reduced lattice basis approach.
The solutions times in \cite{AardalLenstra2004} and our solutions times
are all less than a second, and so it would be interesting to compare the two
methods on larger problems with a significant computation time.

This approach for computing a feasible solution of a fiber
could potentially be used when computing a truncated Markov basis by the
project and lift algorithm since during the algorithm,
we check whether $\f' \in \B{\Lattice}{\f}$ which is the feasibility problem
$\FEASIBLE{\Lattice}{\f-\f'} \ne \emptyset$.
Note that the feasibility approach is well-suited to computing feasibility for
many different fibers simultaneously.
It would be interesting to see the performance of this approach.

\section{Optimality}
\label{Section:Optimality}
In this section, we discuss the use of Gr\"obner bases to solve the integer
program
\[\IP{\Lattice}{c}{\f} := \min \{cx : x \in \FEASIBLE{\Lattice}{\f}\}.\]
The most straight-forward way to solve $\IP{\Lattice}{c}{\f}$ is to compute a
Markov basis of $\Lattice$ and a feasible solution of
$\FEASIBLE{\Lattice}{\f}$ and then compute a $\prec_c$-Gr\"obner basis of
$\Lattice$ for some term order $\prec$, and finally, compute the normal form of
the feasible solution giving the optimal solution.
Here, we are actually solving $\IP{\Lattice}{\prec_c}{\f}$ which
is essentially the same as solving $\IP{\Lattice}{c}{\f}$.
With this method, if we want to solve $\IP{\Lattice}{c}{\f'}$ for a finite
number of $\f'$, we only to compute a feasible solution of
$\FEASIBLE{\Lattice}{\f}$ and redo the normal form computation without needing
to recompute the Markov basis or the Gr\"obner basis.
Also, note that the feasible solutions can be computed at the same time as
computing the Markov basis without much additional computational overhead
(see Section \ref{Section:Feasibility}).

If we wish to solve $\IP{\Lattice}{c}{\f}$ for just one $\f$, then we should
use information specific to that fiber to solve the problem. We can compute a
$\f$-truncated Markov basis $\Lattice$ and a $\f$-truncated
$\prec_c$-Gr\"obner basis of $\Lattice$.
The problem with this method is that we must compute the entire $\f$-truncated
$\prec_c$-Gr\"obner basis in order to prove optimality.
Moreover, if the feasible set $\FEASIBLE{\Lattice}{\f}$ is large, then the
truncated Gr\"obner basis may still be quite large and in some cases as large as
the non-truncated Gr\"obner basis.
Thus, we need further ways of reducing its size.

Some of the non-negativity constraints on the $x$ variables may not strictly be
necessary to solve the problem since they may be redundant or not active near
the optimal solution.
So, we consider relaxations of $\IP{\Lattice}{c}{\f}$ in which we
relax the non-negativity constraints on some of the $x$ variables.
Consider the problem
\[\IPX{\sigma}{\Lattice}{c}{\f} := \min \{cx : x \equiv \f \pmod{\Lattice},
x_\sigmabar \ge 0, x \in \Z^n\}\]
where $\sigma \subseteq \{1,....,n\}$.
Here, we have relaxed the non-negativity constraints on $x_\sigma$.
If $\sigma$ is the set of basic variables given by solving the linear relaxation
of $\IP{\Lattice}{c}{\f}$ using the simplex algorithm (see
\cite{NemhauserWolsey1988}), the relaxation $\IPX{\sigma}{\Lattice}{c}{\f}$ is
called a \emph{group relaxation} (see \cite{Gomory1965}).
Note that, for group relaxations, $\IPX{\sigma}{\Lattice}{c}{\f}$ has an optimal
solution and $\ker(\pi_\sigma) \cap \Lattice = \{\zero\}$.
If $\sigma$ is any subset of the above set for the group relaxation, then
$\IPX{\sigma}{\Lattice}{c}{\f}$ is called an
\emph{extended group relaxation} (see \cite{Wolsey:71}).
Thus, the original problem
$\IP{\Lattice}{c}{\f}=\IPX{\emptyset}{\Lattice}{c}{\f}$ is an
extended group relaxation.

We want to solve these extended group relaxations using Gr\"obner bases and so
we must rewrite $\IPX{\sigma}{\Lattice}{c}{\f}$ in the form
$\IP{\Lattice'}{c'}{\f'}$
for some lattice $\Lattice'$, some right hand side $\f'$, and some cost function
$c'$.
Firstly, any extended group relaxation $\IPX{\sigma}{\Lattice}{c}{\f}$ that has
an optimal solution can always be rewritten in the equivalent form
$\IPX{\sigma}{\Lattice}{\tilde{c}}{\f}$ for some $\tilde{c}$
where $\tilde{c}_\sigma = \zero$, that is, $\tilde{c}_i=0$
for all $i \in \sigma$  (as given by the simplex algorithm
\cite{NemhauserWolsey1988}) where $\IPX{\sigma}{\Lattice}{\tilde{c}}{\f}$ has
the same optimal solution as $\IPX{\sigma}{\Lattice}{c}{\f}$, although the
optimal value may differ by a known constant.
Now consider the projection of
$\IPX{\sigma}{\Lattice}{\tilde{c}}{\f}$ onto the $\sigmabar$ components:
\[ \IP{\Lattice^\sigma}{\tilde{c}_\sigmabar}{\f_\sigmabar} := 
\min \{\tilde{c}_\sigmabar x :x\equiv \f_\sigmabar\pmod{\Lattice^\sigma}, x \in
\N^{|\sigmabar|}\}.\]
Any feasible solution of
$\IP{\Lattice^\sigma}{\tilde{c}_\sigmabar}{\f_\sigmabar}$ lifts to a
feasible solution of $\IPX{\sigma}{\Lattice}{\tilde{c}}{\f}$.
Let $x \in \FEASIBLE{\Lattice^\sigma}{\f_\sigmabar}$, then
$x' := (x-\pi_\sigmabar(\f))+\f \in \FEASIBLEX{\sigma}{\Lattice}{\f}$.
Moreover, an optimal solution
of $\IP{\Lattice^\sigma}{\tilde{c}_\sigmabar}{\f_\sigmabar}$ lifts to an
optimal solution of $\IPX{\sigma}{\Lattice}{\tilde{c}}{\f}$ since
$\tilde{c}_\sigmabar x = \tilde{c} x'$.
Hence, $\IPX{\sigma}{\Lattice}{c}{\f}$, $\IPX{\sigma}{\Lattice}{\tilde{c}}{\f}$
and $\IP{\Lattice^\sigma}{\tilde{c}_\sigmabar}{\f_\sigmabar}$ are all
essentially equivalent problems.

The basic idea of the algorithm is that we start with the group relaxation
$\IP{\Lattice^\sigma}{c_\sigmabar}{\f_\sigmabar}$ where $\sigma$ is defined
as above and we assume $c_\sigmabar = \zero$ without loss of generality.
We then solve the group relaxation
$\IP{\Lattice^\sigma}{c_\sigmabar}{\f_\sigmabar}$.
If the optimal solution of $\IP{\Lattice^\sigma}{c_\sigmabar}{\f_\sigmabar}$
lifts to a feasible solution of $\IP{\Lattice}{c}{\f}$,
then it is optimal for $\IP{\Lattice}{c}{\f}$ and we are done.
Otherwise, we add a non-negativity constraint on one of the unconstrained $x$
variables, that is, we choose $i \in \sigma$ and set
$\sigma := \sigma\setminus\seti$, and solve the extended group relaxation
$\IP{\Lattice^\sigma}{c_\sigmabar}{\f_\sigmabar}$.
Again if the optimal solution of
$\IP{\Lattice^\sigma}{c_\sigmabar}{\f_\sigmabar}$
lifts to a feasible for $\IP{\Lattice}{c}{\f}$, then we are done.
Otherwise, we again add a non-negativity constraint on one of the unconstrained
$x$ variables and solve the new extended group relaxation, and so on, until 
$\IP{\Lattice}{c}{\f}$ is solved. The algorithm must terminate with a solution
because in the worst case we end up solving the original problem
$\IP{\Lattice}{c}{\f}$ ($\sigma=\emptyset$).

To solve an extended group relaxation
$\IP{\Lattice^\sigma}{c_\sigmabar}{\f_\sigmabar}$, we first compute a
$\f_\sigmabar$-truncated Markov basis of $\Lattice^\sigma$ and a feasible
solution of $\FEASIBLE{\Lattice^\sigma}{\f_\sigmabar}$.
Secondly, we compute a $\f_\sigmabar$-truncated
$\prec^\sigma_{c_\sigmabar}$-Gr\"obner basis of $\Lattice^\sigma$.
Then, we compute the normal form of the feasible solution giving the optimal
solution of $\FEASIBLE{\Lattice^\sigma}{\f_\sigmabar}$.

Initially, we need to compute a $\f_\sigmabar$-truncated Markov basis of
$\Lattice^\sigma$ for the group problem.
Let $B$ be a basis of the lattice $\Lattice$.
Then for the group problem, $\sigma$ gives $\rank(B)$ linear independent columns
of $B$.
Thus, as discussed at the end of Section 
\ref{Section: Truncated Markov bases}, we can compute a Markov basis
via a HNF computation.
Similarly, we can compute a feasible solution as discussed at the end of Section
\ref{Section:Feasibility}.

At each iteration of the algorithm, we could compute a
$\f_\sigmabar$-truncated Markov basis of $\Lattice^\sigma$ starting from
scratch each time as described in Section
\ref{Section: Truncated Markov bases},
but instead, we can compute it incrementally exactly as in the
Project-and-Lift algorithm.
In the previous iteration, we will have computed a $\f_{\sigmabar'}$-truncated
Markov basis of $\Lattice^{\sigmabar'}$ where $\sigma' := \sigma\cup\seti$ for
some $i \in \sigmabar$.
Hence, using Lemmas \ref{Lemma:project-and-lift} and
\ref{Lemma: gen non-negative vector}, we can compute a
$\f_\sigmabar$-truncated Markov basis of $\Lattice^\sigma$ from a
$\f_{\sigmabar'}$-truncated Markov basis of $\Lattice^{\sigmabar'}$ in one
step. In effect, we perform the Project-and-Lift algorithm
simultaneously.
Also, as discussed in Section \ref{Section:Feasibility}, we can compute a
feasible solution of $\IP{\Lattice^\sigma}{c_\sigmabar}{\f_\sigmabar}$ from a
feasible solution of $\IP{\Lattice^{\sigma'}}{c_{\sigmabar'}}{\f_{\sigmabar'}}$.

See Algorithm \ref{Algorithm: IP} for a description of the optimisation
algorithm.

\begin{algorithm}[ht]
\caption{Optimisation algorithm}
\begin{algorithmic}
\REQUIRE an integer program $\IP{\Lattice}{c}{\f}$.
\ENSURE an optimal solution $x \in \FEASIBLE{\Lattice}{\f}$ or
  \emph{infeasible}.
\STATE Find a set $\sigma\subseteq\{1,\ldots,n\}$ such that
$\ker(\pi_\sigma) \cap \Lattice = \{\zero\}$.
\WHILE{$\sigma\neq\emptyset$}
  \STATE Compute a set $M \subseteq \Lattice^\sigma$ that is a
    $\f_\sigmabar$-truncated Markov basis of $\Lattice^\sigma$.
  \STATE Compute a feasible solution
    $x\in\FEASIBLE{\Lattice^\sigma}{\f_\sigmabar}$ and \textbf{return}
    \emph{infeasible} if $\FEASIBLE{\Lattice^\sigma}{\f_\sigmabar}=\emptyset$.
  \STATE $G:= \CP_{\f_\sigmabar}(\prec^\sigma_{c_\sigmabar},M)$
  \STATE $x:= \NF(x,G)$
  \STATE $x:= \pi_\sigma^\inv(x - \f_\sigmabar) + \f$
  \STATE \Kw{if} $x \ge 0$ \Kw{then} \Kw{return} $x$ 
  \STATE Select $i \in \sigma$.
  \STATE $\sigma:=\sigma\setminus\seti$
\ENDWHILE
\end{algorithmic}
\label{Algorithm: IP}
\end{algorithm}

At each iteration, we must select the next $i$. An obvious choice is to
select the component with the most violated non-negativity constraint, that is,
the most negative component.

Existing branch-and-bound methods for integer programming can take
advantage of a good feasible solution of $\IP{\Lattice}{c}{\f}$, but the above
method cannot.
However, we can take advantage of a good feasible solution
since a feasible solution gives us an upper bound on $\IP{\Lattice}{c}{\f}$ that
can be used to strengthen truncation.

Consider the following reformulation of $\IP{\Lattice}{c}{\f}$ using some upper
bound $k\in\Z$ on $\IP{\Lattice}{c}{\f}$:
\begin{align*}
\IP{\Lattice}{c}{\f}
:=& \min \{cx : x \equiv \f \pmod{\Lattice}, x\in\N^n\} \\
 =& \min \{cx : x \equiv \f \pmod{\Lattice}, cx \le k, x\in\N^n\} \\
 =& \min \{-y : x \equiv \f \pmod{\Lattice}, cx+y=k, x\in\N^n, y\in\N\} + k.
\end{align*}
To solve the reformulation, we first need to express it in the form
$\IP{\Lattice'}{c'}{\f'}$ for some lattice $\Lattice'$, right hand side $\f'$
and some cost function $c'$.
Let $\Lattice' := \{ (u,-cu) : u \in \Lattice\}$, $\f':=(\f,k-c\f)$, and
$c' = (\zero,-1)$.
Then,
\begin{align*}
\IP{\Lattice'}{c'}{\f'}
:=&\min\{-y:(x,y) \equiv (\f,k-c\f) \pmod{\Lattice'},(x,y)\in\N^{n+1}\}\\
=&\min\{-y:(x-\f,y-(k-c\f)) \in \Lattice',(x,y)\in\N^{n+1}\}\\
=&\min\{-y:(x-\f,y-(k-c\f))=(u,-cu), u \in \Lattice,(x,y)\in\N^{n+1}\}\\
=&\min\{-y:x\equiv \f\pmod{\Lattice}, cx+y=k, x\in\N^n,y\in\N\}.
\end{align*}
Thus, $\IP{\Lattice}{c}{\f} = \IP{\Lattice'}{c'}{\f'}+k$.
So, we can solve $\IP{\Lattice}{c}{\f}$ by solving
$\IP{\Lattice'}{c'}{\f'}$ using the methods discussed previously.

Hopefully, a $\f'$-truncated $\prec_{c'}$-Gr\"obner basis of $\Lattice'$ is a
lot smaller than a $\f$-truncated Gr\"obner basis of $\Lattice$.
Crucially, the size of a minimal $\f'$-truncated $\prec_{c'}$-Gr\"obner basis of
$\Lattice'$ cannot exceed the size of a minimal $\f$-truncated
$\prec_c$-Gr\"obner basis, and so, we are not computing more than before.
This follows from Lemma \ref{Lemma:project-and-lift}; since
$\pi_{n+1}(\Lattice') = \Lattice$, $\pi_{n+1}(\f') = \f$,
$c'=\ebar^{n+1}$, and $c = \omegabar^{n+1}$,
if $G$ is a $\f$-truncated
$\prec_c$-Gr\"obner basis of $\Lattice$, then $\pi^\inv_{n+1}(G)$
is a $\f'$-truncated $\prec_{c'}$-Gr\"obner basis of $\Lattice'$.
How much smaller a $\f'$-truncated $\prec_{c'}$-Gr\"obner basis of $\Lattice'$
is than a $\f$-truncated Gr\"obner basis of $\Lattice$
will depend on the strength of the upper bound.
It is potentially just the empty set.

If we are given a feasible solution $x \in \FEASIBLE{\Lattice}{\f}$,
then we can set $k = cx-1$.
If $x$ is the optimal solution of $\IP{\Lattice}{c}{\f}$,
then setting $k=cx-1$, we have $\FEASIBLE{\Lattice'}{\f'} = \emptyset$.
Therefore, a minimal $\f'$-truncated Markov basis of $\Lattice'$ is empty and
a minimal $\f'$-truncated $\prec_{c'}$-Gr\"obner basis of $\Lattice'$ is also
empty.
So, potentially, computing a $\f'$-truncated $\prec_{c'}$-Gr\"obner basis of
$\Lattice'$ is a lot more efficient than computing a $\f$-truncated
$\prec_c$-Gr\"obner basis of $\Lattice$;
we only have to compute an empty set to show optimality!

Even if the bound $k$ is not very good and thus does not help
truncation much, it is still definitely worthwhile
to solve $\IP{\Lattice'}{c'}{\f'}$ instead of $\IP{\Lattice}{c}{\f}$.
The reason is that by introducing the constraint $cx \le k$ into the problem,
more components may become \emph{bounded} and thus the Gr\"obner basis and
the Markov basis computations for $\IP{\Lattice'}{c'}{\f'}$ 
are faster (see the section in \cite{Hemmecke+Malkin:2006} on Criterion 2).
The computations are also faster for the extended group relaxations of
$\IP{\Lattice'}{c'}{\f'}$ as well.

\begin{example}
Let $\Lattice = \Lattice_A$ for the matrix $A$ given in Example
\ref{Example: GB}.
In the following table, we list the time taken to compute the optimal solution
of $\IP{\Lattice}{c}{\f}$ for different $\f$'s given a feasible solution where
$c=(3,15,1,5,2,17,16,16,15,9,7,11,13)$.
\begin{center}
\begin{tabular}{|c|c|c|c|}
\hline
$\f$ & Group Relaxation & Final Relaxation & Time \\
\hline
\hline
$(1,1,1,0,1,0,1,1,0,1,0,1,0)$& $\{3,4,5,10\}$& $\{3,4,5,10\}$ & $0.28$s \\
\hline
$(1,0,1,0,3,0,1,5,0,1,0,9,0)$& $\{1,5,8,12\}$& $\{1,5,8,12\}$ & $0.02$s \\
\hline
$(1,1,1,1,1,1,1,1,1,1,1,1,1)$& $\{3,4,5,10\}$& $\{5,10\}$   & $0.61$s \\
\hline
$(1,2,0,3,5,0,1,3,0,4,0,1,0)$& $\{3,4,5,10\}$& $\{3,5,10\}$   & $0.46$s \\
\hline
$(19,7,3,8,13,11,1,15,4,8,17,9,5)$& $\{3,4,5,10\}$ & $\{3,4,5,10\}$ & $0.28$s\\
\hline
\end{tabular}
\end{center}
Here, we used the quick check for truncation.
\end{example}

%

If we do not have a good feasible solution or any feasible solution at all
available, we can still use the extended formulation.
Assume that we are given a lower bound $l$ on the optimal value, which we can
always find by solving the linear relaxation.
We then try to find a feasible solution by computing a $\f'$-truncated
Markov basis of $\Lattice'$ where $k=l$.
If we find a feasible solution, then it must be
optimal. Otherwise, we recompute a $\f'$-truncated
Markov basis of $\Lattice'$ where $k=l+1$ and again try to find a feasible
solution. We repeat this procedure by incrementing $k$ until we find a feasible
solution which must be an optimal solution.
This procedure has the advantage that we only compute feasible solutions and not
optimal solutions, and thus, we avoid some Gr\"obner basis computations.


The recomputation of $\f'$-truncated Markov bases to find initial solutions
initially seems inefficient; however, this is not the case because we can reuse
the previous computations.
Let $\f' = (\f,k-c\f)$ and $\f'' = (\f,k+1-c\f)$.
Then $\f' \in \B{\Lattice}{\f''}$ since
$(\zero,1) \in \FEASIBLE{\Lattice'}{\f''-\f'}$ and so
$\B{\Lattice}{\f'} \subseteq \B{\Lattice}{\f''}$.
Hence, to compute a $\f''$-truncated Gr\"obner basis or Markov basis 
of $\Lattice$ requires also computing a $\f'$-truncated Gr\"obner basis or
Markov basis respectively anyway.
This applies not only at the final stage of the algorithm for each value of $k$,
but also at each intermediate stage for the extended group relaxations of
$\IP{\Lattice'}{c'}{\f'}$.
So, it requires keeping all the intermediate stages
around, and is thus more
complex than the first algorithm presented in this section.

We have not implemented such an approach yet.  It would be interesting to see
how it performs. However, if we know a very good initial feasible solution, we
would expect that the previous method is faster.

\section{Conclusion}
We have demonstrated that it is possible to
significantly improve upon the
performance of previous Gr\"obner basis based approaches.
However, Gr\"obner basis approaches have not yet proven to
be competitive with traditional branch-and-bound based approaches to 
integer programming for industrial applications, but given the
significant advances shown here, perhaps with further research,
Gr\"obner basis methods will be useful for some classes of problems
when combined with existing methods.

\appendix
\section{Lattice Programs and Integer Programs}
\label{Appendix: Lattice Programs and Integer Programs}
In this appendix, we show the equivalence of fibers and feasible
sets of integer programs and of lattice programs and integer programs.

Consider the set
$\{x:Ax=b, x_\sigmabar \geq \zero, x\in\Z^n\}$
where $A \in \Z^{m\times n}$, $b \in \Z^n$, $\sigma \subseteq \{1,...,n\}$,
$\sigmabar$ is the complement of $\sigma$, and $x_\sigmabar$ is the set of
variables indexed by $\sigmabar$.
Thus, the $x_\sigmabar$ variables are non-negative and the $x_\sigma$ variables
(the $x$ variables indexed by $\sigma$) are unrestricted in sign.
The set of feasible solutions for any integer program can be represented in
this form.
We will rewrite this set as a lattice fiber.
We will actually only consider the projection of this set onto the non-negative
variables (the $\sigmabar$ components), that is,
$F^\sigma_A(b) := \{x_\sigmabar:Ax=b, x_\sigmabar \geq \zero, x\in\Z^n\}$
since any solution of $F^\sigma_A(b)$ can easily be extended to a solution in
the original space;
more specifically, given $\tilde{x}_\sigmabar \in F^\sigma_A(b)$,
we can find
$\tilde{x}_\sigma\in\Z^{|\sigma|}$ where
$A_\sigma \tilde{x}_\sigma=b-A_\sigmabar \tilde{x}_\sigmabar$ by
using the HNF algorithm
where $A_\sigma$ and $A_\sigmabar$ are the sub-matrices of $A$ whose
columns are indexed by $\sigma$ and $\sigmabar$ respectively.

We now write $F^\sigma_A(b)$ as a fiber of a lattice.
Let $\Lattice_A = \{u \in \Z^n: Au = \zero\}$,
and let $\f \in \Z^n$ where $A\f=b$.
Note that we can compute a basis of $\Lattice_A$ and a $\f$ using the 
HNF algorithm.
Also, let $\Lattice^\sigma_A$ be the
lattice $\Lattice_A$ projected onto the $\sigmabar$ components
(i.e. $\Lattice^\sigma_A := \{u_\sigmabar:Au=\zero\}$).
Note that a projection of a lattice is always again a lattice.
It follows that
\begin{align*}
F^\sigma_A(b)   &= \{x_\sigmabar :Ax=b, x_\sigmabar \geq \zero, x\in\Z^n\}\\
                &= \{x_\sigmabar :Ax=A\f, x_\sigmabar\geq \zero, x\in\Z^n\}\\
                &= \{x_\sigmabar :A(x-\f)=\zero, 
                x_\sigmabar\geq \zero, x\in\Z^n\}\\
                &= \{x_\sigmabar :x-\f \in \Lattice_A, x_\sigmabar \geq \zero,
                        x\in\Z^n\}\\
                &= \{x_\sigmabar :x_\sigmabar -\f_\sigmabar\in\Lattice^\sigma_A,
                     x_\sigmabar \geq \zero, x_\sigmabar \in\Z^{|\sigmabar|}\}\\
                &= \FEASIBLE{\Lattice^\sigma_A}{\f_\sigmabar}.
\end{align*}
So, the set of feasible solutions to any integer program can be represented as a
fiber of a lattice.

Conversely, given a $\Lattice$ and a $\f\in\Z^n$, we can represent the fiber
$\FEASIBLE{\Lattice}{\f}$ as the set of feasible solutions to an integer
program.
Let $S\in\Z^{n \times k}$ be a matrix where the columns of $S$ span the lattice
$\Lattice$. Then,
\begin{align*}
\FEASIBLE{\Lattice}{\f}&=\{x\in\N^n: x\equiv \f\pmod{\Lattice}\}\\
                      &=\{x\in\N^n: x-\f \in \Lattice\}\\
                      &=\{x\in\N^n: x-\f = Sy, y \in \Z^k\}\\
                      &=\{x: (I,-S)(x,y)=\f,  x\ge0, x\in\Z^n,y \in \Z^k\}.\\
                      &=\{x: A(x,y)=\f,  x\ge0, x\in\Z^n,y \in \Z^k\}.\\
                      &= F^\sigma_A(\f)
\end{align*}
where $A=(I,-S)\in\Z^{(n+k)\times n}$ and $\sigma =\{1,...,n\}$ and
$\sigmabar=\{n+1,...,n+k\}$.
Thus, fibers of lattices and feasible sets of integer programs are two different
representations of the same set and so are equivalent concepts.

Next we show the equivalence of lattice programs and integer programs.
Any integer linear program can be written in the form
$\min\{cx:Ax=b, x_\sigmabar \geq \zero, x\in\Z^n\}$
where $A \in \Z^{m\times n}$, $b\in\Z^m$, and $c\in\Z^n$.
If this integer program has an optimal solution, then there exists a
$\tilde{c} \in \Z^n$ where  $\tilde{c}_\sigma = \zero$ such that
$cx=\tilde{c}x+k$ for some constant $k\in\Q^n$ and every feasible solution $x$. 
Then, 
\begin{align*}
\min\{cx:Ax=b, x_\sigmabar \geq \zero, x\in\Z^n\} &= 
\min\{\tilde{c}x+k:Ax=b, x_\sigmabar \geq \zero, x\in\Z^n\} \\
&=\min\{\tilde{c}_\sigmabar x_\sigmabar:Ax=b, x_\sigmabar \geq \zero, x\in\Z^n\}+k\\
&=\min\{\tilde{c}_\sigmabar x_\sigmabar:x_\sigmabar \in F^\sigma_A(b)\}+k\\
&=\min\{\tilde{c}_\sigmabar x_\sigmabar:
x_\sigmabar \in \FEASIBLE{\Lattice^\sigma_A}{\f_\sigmabar}\}+k
\end{align*}
where $\Lattice^\sigma_A$ and $\f_\sigmabar$ are defined as above.
More specifically,
we can always find a $\tilde{c} \in \Z^n$ where  $\tilde{c}_\sigma = \zero$ such
that $cx=\tilde{c}x+k$ for every feasible solution $x$
if and only if $cx = 0$ for all $x \in \Z^n$ where $x_\sigmabar = \zero$ and
$Ax=\zero$. Note that if $cx \ne 0$ for some $x \in \Z^n$ where
$x_\sigmabar=\zero$
and $Ax=\zero$, then the integer program has no optimal solution.
We can check this condition and find a valid $\tilde{c}$ easily using the
HNF algorithm.

Hence, any integer program with an optimal solution can be written in the form
\[\IP{\Lattice}{c}{\f}:=\min\{cx:x\in\FEASIBLE{\Lattice}{\f}\}\]
for some lattice $\Lattice \subseteq \Z^n$, a vector $\f \in \Z^n$, and a vector
$c \in \Z^n$.
As we saw at the end of Section \ref{Section: Truncated Groebner bases},
solving $\IP{\Lattice}{c}{\f}$ is equivalent to
solving the \emph{lattice program} $\IP{\Lattice}{\succ}{\f}$
where $\succ$ is compatible with $c$

\section*{Acknowledgments}
I would like to thank Raymond Hemmecke and Laurence Wolsey
for many fruitful discussions and helpful comments. 

\bibliographystyle{plain}

\end{document}